\documentclass[12pt]{amsart}
\usepackage{amssymb}
\usepackage{amsmath}
\usepackage{mathabx}
\usepackage{amsthm}
\usepackage{amstext}
\usepackage{amsopn}
\usepackage{mathrsfs}
\usepackage{latexsym}
\DeclareMathAlphabet{\mathpzc}{OT1}{pzc}{m}{it}
\allowdisplaybreaks
\newtheorem{Definition}{Definition}[section]
\newtheorem{Proposition}{Proposition}[section]
\newtheorem{Lemma}{Lemma}[section]
\newtheorem{Theorem}{Theorem}[section]
\newtheorem{Corollary}{Corollary}[section]
\newtheorem{Remark}{Remark}[section]
\newtheorem{Example}{Example}[section]

\begin{document}
\bibliographystyle{plain}
\footnotetext{
This work is supported by Narodowe Centrum Nauki grant No. 2014/15/B/ST1/00166}
\title[Matricial circular systems and random matrices]
{Matricial circular systems and random matrices}
\author[R. Lenczewski]{Romuald Lenczewski}
\address{Romuald Lenczewski, \newline
Wydzia\l{} Matematyki, Politechnika Wroc\l{}awska, \newline
Wybrze\.{z}e Wyspia\'{n}skiego 27, 50-370 Wroc{\l}aw, Poland  \vspace{10pt}}
\email{Romuald.Lenczewski@pwr.edu.pl}
\begin{abstract}
We introduce and study {\it matricial circular systems} of operators which 
play the role of matricial counterparts of circular operators. 
They describe the asym\-ptotic joint *-distributions of blocks of independent block-identically 
distributed Gaussian random matrices with respect to partial traces.
Using these operators, we introduce {\it circular free Meixner distributions} as the non-Hermitian counterparts 
of free Meixner distributions and construct for them a random matrix model.
Our approach is based on the concept of matricial freeness applied to operators 
on Hilbert spaces. It is closely related to freeness with amalgamation over the algebra $A$ of $r\times r$ diagonal matrices applied to operators on Hilbert $A$-bimodules.
\end{abstract}

\maketitle
\section{Introduction}
Circular and semicircular systems of operators were introduced by Voiculescu [16] in 
the context of random matrices and their asymptotics under the expectation of 
the normalized trace and applied to isomorphisms of free group factors. 
This result followed his fundamental asymptotic freeness result for independent 
Hermitian and non-Hermitian Gaussian random matrices with complex i.i.d. entries, respectively [15]. Decomposition of these matrices into blocks corresponds to decompositions 
of circular and semicircular operators into circular and semicircular systems,
respectively. 

In order to describe asymptotic joint *-distributions of Gaussian random matrices 
with non-identically distributed entries, one has to use a more general framework.
The first approach was given by Shlyakhtenko [12], who used free probability with operator-valued 
states and the associated notion of freeness with amalgamation
over the algebra of diagonal matrices to describe them. 
This approach was further developed by Benaych-Georges 
[1] in his study of asymptotic joint *-distributions of blocks of 
Gaussian random matrices, including the case when one of the asymptotic dimensions vanishes.

Recently, we proposed a different approach to study the case when 
the considered independent Gaussian (in the Hermitian case, 
not only Gaussian) random matrices have independent block-identically
distributed complex entries, which we abbreviate {\it i.b.i.d.} [6,7]. We showed that
in order to describe the asymptotic joint (*-) distributions 
of symmetric blocks of these matrices under {\it partial traces} 
(understood as expectations of normalized traces over subsets of indices) 
rather than under the expectation of the normalized complete trace, 
one needs to replace (*-) free semicircular or circular
operators by their matricial counterparts, related to each other by 
a matricial generalization of (*-) freeness with respect to an array of scalar-valued states. 

The underlying concept of independence is that of {\it matricial freeness} [5], 
which can be described by means of the intuitive equation
$$
matricial\;freeness=freeness \;\;\&\;\;matriciality,
$$
and its symmetrized version called {\it symmetric matricial freeness} [6] which 
also played an important role in our previous developments. 
A connection between these concepts and blocks of large random matrices was first 
given for one matrix [6] and then for families of independent matrices [7].

Using this concept, we study the asymptotic joint *-distributions
of blocks of independent Gaussian random matrices under partial traces.
If we are given an ensemble of independent non-Hermitian $n\times n$ random matrices
$\{Y(u,n):u\in \mathpzc{U}\}$ whose entries are suitably normalized i.b.i.d. complex 
Gaussian random variables with zero mean for each natural $n$, then 
we show that the mixed *-moments of their blocks $\{S_{p,q}(u,n):1\leq p, q \leq r, u\in \mathpzc{U}\}$ 
converge under partial traces $\tau_q(n)$ to the mixed *-moments of 
certain bounded non-self-adjoint operators under 
scalar-valued states $\Psi_{q}$, which we write informally
$$
\lim_{n\rightarrow \infty}S_{p,q}(u,n)=\zeta_{p,q}(u),
$$
and the operators $\zeta_{p,q}(u)$ are called {\it matricial circular operators}. The whole family
$\{\zeta_{p,q}(u):1\leq p, q \leq r, u\in \mathpzc{U}\}$ will be called a {\it matricial circular system}. 
Let us remark that our approach allows us to treat 
all rectangular blocks, including those which are {\it unbalanced} or {\it evanescent} [7].

Matricial circular operators and their mixed *-moments under the states $\Psi_q$ 
are the main objects of our study. Their main properties are listed below.
\begin{enumerate}
\item
There is a connection between our approach and operator-valued free probability [17].
In particular, the operators $\zeta_{p,q}(u)$ can be identified with 
elements of the Toeplitz algebra ${\mathcal T}_{{\mathcal H}}$ 
of the form
$$
g_{p,q}(u)=F_pg(u)F_q,
$$
where ${\mathcal H}$ is a Hilbert $A$-bimodule with $A$ being the $C^*$-algebra 
of $r\times r$ diagonal matrices with canonical generators $F_1, \ldots , F_r$,
and $\{g(u):u\in \mathpzc{U}\}$ is the family of *-free $A$-valued circular elements 
under the canonical conditional expectation $E:{\mathcal T}_{{\mathcal H}}\rightarrow A$ 
used by Shlyakhtenko [12, Theorem 4.1] in his $A$-valued approach to random band matrices.
\item
The joint *-distributions of the operators $\zeta_{p,q}(u)$ under the
states $\Psi_q$ can be identified with the joint *-distributions of elements 
$g_{p,q}(u)$ under the compressed conditional expectation 
$E_q:{\mathcal T}_{{\mathcal H}}\rightarrow A$ given by 
$E_q(x)=E(F_qxF_q)$. More generally, their joint *-distribution under the convex
linear combination $\Psi=\sum_{q=1}^{r}d_{q}\Psi_q$ can be written in the form
$$
\Psi(\widetilde{\sigma}(x))={\rm Tr}(E(x)D)
$$
where $D={\rm diag}(d_1, \ldots, d_r)$, $x$ is a word in elements $g_{p,q}(u)$ 
and $\widetilde{\sigma}$ is a unital *-homomorphism which maps $g_{p,q}(u)$ onto
$\zeta_{p,q}(u)$.
\item
The symmetrizations
$$
\eta_{p,q}(u)=\left\{\begin{array}{ll}\zeta_{p,q}(u) +\zeta_{q,p}(u) & {\rm if}\; p\neq q\\
\zeta_{q,q}(u) & {\rm if}\;p=q
\end{array}
\right.,
$$
give the operatorial realization of asymptotic symmetric blocks [7].
We will call them {\it symmetrized matricial circular operators}, in contrast to [7], where 
we used the name `matricial circular operators'.
Each $\eta_{p,q}(u)$ has the standard 
circular distribution under both $\Psi_q$ and $\Psi_p$,
provided the covariances of the summands are equal to one.
\end{enumerate}

In general, one of our main goals in [6,7] and in the present paper is to construct a unified system of Hilbert space operators describing the asymptotic *-distributions of a large class of random matrices.
One of the benefits is the construction of new random matrix models, like that for
free Meixner laws [8], the idea which is continued in this paper. Namely, using matricial circular systems, we define operators whose distributions play the role of non-Hermitian counterparts of free Meixner distributions called {\it circular free Meixner distributions} and construct for them a random matrix model. Some implications concerning the triangular random matrix models studied by Dykema and Haagerup [3,4] will be given elsewhere.  

The paper is organized as follows. In Section 2, we recall the definition of matricially free creation operators and we show a connection between these operators and 
creation operators in certain Hilbert $A$-bimodules, where $A$ is the $C^*$-algebra of $r\times r$ diagonal matrices.  
In Section 3, we introduce the matricial circular operators and establish their
realization as operator-valued matrices.
In Section 4, we show that *-distributions of blocks of i.b.i.d. GRM 
under partial traces converge to the *-distributions of the matricial circular operators.
In Sections 5-7, we study the mixed *-moments and mixed *-cumulants of 
these operators, as well as the corresponding moment series and cyclic $R$-transform.
In Section 8, we introduce circular free Meixner distributions, for which we construct 
a random matrix model.

When speaking of only one matrix and one array of operators for some fixed $u$, we will omit $u$ in our notations, writing $Y(n)$, $S_{p,q}(n)$, $(\zeta_{p,q})$, etc.

\section
{Matricial freeness and Hilbert $A$-bimodules}

Our study of the asymptotic distributions of random matrices in [6,7] 
was based on the construction of the appropriate Hibert space of Fock type,
in which `matriciality' is added to `freeness' and in which certain 
systems of operators, called matricially free Gaussian operators, replace semicircular systems.  
In this section, we recall these notions and show their relation to 
the free Fock space over the Hilbert $A$-bimodule and $A$-valued semicircular systems, 
where $A$ is the $C^*$-algebra of $r\times r$ complex diagonal matrices.

Let us recall the Hilbert space setting. Let 
$\mathpzc{J}\subseteq [r]\times [r]$, where 
$[r]:=\{1,2, \ldots, r\}$, although we will mainly deal with the situation 
when these sets are equal. However, passing from $[r]\times [r]$ to any proper subset presents no difficulty
and can be achieved by setting certain operators to be zero. Moreover, we will consider another finite set of 
indices $\mathpzc{U}$, setting for convenience $\mathpzc{U}=[t]$ for some integer $t$.
To each $(p,q)\in \mathpzc{J}$ and $u\in \mathpzc{U}$ we then associate a Hilbert space 
${\mathcal H}_{p,q}(u)$. Using this family of Hilbert spaces, we can construct
our Fock space, a matricial version of the free Fock space. 

\begin{Definition}
{\rm By the {\it matricially free Fock space of tracial type} over the family 
of Hilbert spaces $\{{\mathcal H}_{p,q}(u):(p,q)\in \mathpzc{J}, u\in \mathpzc{U}\}$ 
we understand the Hilbert space direct sum
\begin{equation*}
{\mathcal M}= \bigoplus_{q=1}^{r} {\mathcal M}_{q},
\end{equation*}
where each summand is of the form
\begin{equation*}
{\mathcal M}_{q}={\mathbb C}\Omega_{q}\oplus \bigoplus_{m=1}^{\infty}
\bigoplus_{\stackrel{p_1,\ldots , p_m}
{\scriptscriptstyle u_1, \ldots , u_n}
}
{\mathcal H}_{p_1,p_2}(u_{1})\otimes {\mathcal H}_{p_2,p_3}(u_2)
\otimes \ldots \otimes 
{\mathcal H}_{p_m,q}(u_{m}),
\end{equation*}
with $\Omega_q$ being a unit vector, endowed with the canonical inner products.}
\end{Definition}

Observe that the neighboring Hilbert spaces 
can coincide if and only if they have diagonal indices. 
This shows that there is an essential difference between the diagonal Hilbert 
spaces and the off-diagonal ones in this structure. 
 
\begin{Proposition}
If ${\mathcal H}_{p,q}(u)={\mathbb C}e_{p,q}(u)$ for any $p,q,u$, where $e_{p,q}(u)$ is a unit vector, the 
canonical orthonormal basis ${\mathcal B}$ of the matricially free Fock space ${\mathcal M}$ consists of
$$
e_{p_1,p_2}(u_1)\otimes e_{p_2,p_3}(u_2)\otimes  \ldots \otimes e_{p_m,q}(u_{m}),
$$
where $p_1, \ldots, p_m,q\in [r]$, $u_{1},\ldots , u_m\in \mathpzc{U}$ and $m\in {\mathbb N}$, and of vectors $\Omega_1, \ldots , \Omega_r$.
\end{Proposition}
{\it Proof.}
This fact is obvious.
\hfill $\blacksquare$\\

\begin{Definition}
{\rm Let $B(u)=(b_{p,q}(u))$ be an $r\times r$ array of positive real numbers for any $u\in \mathpzc{U}$.
We associate with each such matrix the {\it matricially free creation operators} whose 
non-trivial action onto the basis vectors is  
\begin{eqnarray*}
\wp_{p,q}(u)\Omega_q&=&\sqrt{b_{p,q}(u)}e_{p,q}(u)\\
\wp_{p,q}(u)(e_{q,t}(s))&=&\sqrt{b_{p,q}(u)}(e_{p,q}(u)\otimes e_{q,t}(s))\\
\wp_{p,q}(u)(e_{q,t}(s)\otimes w)&=&\sqrt{b_{p,q}(u)}(e_{p,q}(u)\otimes e_{q,t}(s)\otimes w)
\end{eqnarray*}
for any $p,q,t\in [r]$ and $u,s\in \mathpzc{U}$, where $e_{q,t}(s) \otimes w$ is assumed to be 
a basis vector. Their actions onto the remaining basis vectors give zero. 
The corresponding {\it matricially free annihilation operators} are their adjoints 
denoted $\wp_{p,q}^{*}(u)$. If $b_{p,q}(u)=1$, we will call the associated operators {\it standard}.}
\end{Definition}

Now, let us turn to the Hilbert $A$-bimodule setting, where $A$ is the algebra 
of complex diagonal $r\times r$ matrices. For this purpose, we specify the main objects.
\begin{enumerate}
\item
For the family of Hilbert spaces
$$
\{\mathcal{H}_{p,q}(u): (p,q)\in \mathpzc{J}, u\in \mathpzc{U}\}
$$ 
we form their direct sum
$$
\mathcal{H}=\bigoplus_{u\in \mathpzc{U}}\bigoplus_{(p,q)\in \mathpzc{J}}
{\mathcal H}_{p,q}(u),
$$
and make it into a Hilbert $A,A$-bimodule (or, Hilbert $A$-bimodule), where $A$ is the algebra of diagonal $r\times r$ complex matrices (it is a finite dimensional commutative $C^*$-algebra).
\item
This structure is obtained by the right and left actions of $A$ onto $\mathcal{H}$ given by
$$
ha=(h(u)a)_{u\in \mathpzc U}
$$
$$
ah=(ah(u))_{u\in \mathpzc U}
$$
for $h=(h(u))_{u\in \mathpzc U}\in {\mathcal H}$, respectively, 
where
$$
h(u)a=(h_{p,q}(u)a_{q})_{(p,q)\in \mathpzc J}
$$
$$
ah(u)=(a_{p}h_{p,q}(u))_{(p,q)\in \mathpzc J}
$$
and $a={\rm diag}(a_1, \ldots, a_{r})\in A$. Note that if we treat
$h(u)=(h_{p,q}(u)))_{(p,q)\in\mathpzc J}$ as a matrix, these actions can be viewed as right and left multiplication by $a$.
\item
The $A$-valued inner product on ${\mathcal H}$ is defined 
in terms of ${\mathbb C}$-valued inner products on the considered family of Hilbert spaces 
by the formula
$$
\langle h, g\rangle_{A}=\sum_{u\in \mathcal U}
\sum_{(p,q)\in \mathpzc{J}}\langle h_{p,q}(u), g_{p,q}(u)\rangle F_{q},
$$
where $F_{q}$ is the diagonal $r\times r$ matrix of rank one defined by 
$(F_{q})_{i,j}=\delta_{i,q}\delta_{j,q}$. In particular, the left action 
of $A$ onto ${\mathcal H}$ is a *-homomorphism from $A$ into the $C^*$-algebra ${\mathcal L}({\mathcal H})$ of adjointable linear maps on ${\mathcal H}$, and the right action of $A$ onto ${\mathcal H}$ 
is an antimultiplicative linear map. With these actions and the inner product defined above ${\mathcal H}$ becomes a Hilbert $A$-bimodule. 
\item
Using the Hilbert A-bimodule structure on ${\mathcal H}$, we define the full Fock space over ${\mathcal H}$ as the direct sum
$$
{\mathcal F}({\mathcal H})=\bigoplus_{n=0}^{\infty}{\mathcal H}^{\otimes_{A}n}
$$
where ${\mathcal H}^{\otimes_{A}0}=A$ and 
${\mathcal H}^{\otimes_{A}n}$ is the $n$-fold tensor product 
${\mathcal H}\otimes_{A}\ldots \otimes_{A}{\mathcal H}$. This direct sum is equipped with the canonical $A$-valued inner product, namely: on $A$ we set $\langle F_q, F_s\rangle_{A}=\delta_{q,s}F_{q}$, on ${\mathcal H}^{\otimes_{A}1}$ it is inherited from ${\mathcal H}$, whereas on higher order tensor products it is given by the recursive formula
$$
\langle h_1\otimes_{A}\ldots \otimes_{A}h_{n},
g_1\otimes_{A}\ldots \otimes_{A}g_{n}\rangle_{A}
$$
$$
=\langle h_{n},\langle h_1\otimes_{A}\ldots \otimes_{A}h_{n-1},
g_1\otimes_{A}\ldots \otimes_{A}g_{n-1}\rangle_{A}g_{n}\rangle_{A}
$$
for any $n\geq 1$ and $h_1, \ldots , h_n, g_1, \ldots , g_n\in \mathcal{H}$. Then
${\mathcal F}:={\mathcal F}({\mathcal H})$ 
becomes a Hilbert $A$-bimodule with the natural left and 
right actions of $A$.
\item
For each $h\in {\mathcal H}$, one defines the associated left creation operator $\ell(h)$ by setting 
\begin{eqnarray*}
\ell(h)a&=&ha\\
\ell(h)h_1\otimes_A \ldots \otimes_A h_n&=&h\otimes_Ah_1\otimes_A \ldots \otimes_{A}h_{n} 
\end{eqnarray*}
for any $a\in A$ and any $h_1, \ldots , h_n\in \mathcal{H}$. Then $\ell(h)$ is adjointable on ${\mathcal F}(\mathcal{H})$ and its adjoint is denoted $\ell(h)^*$. 
\item
The following relations hold:
\begin{eqnarray*}
\ell(h)^*\ell(g)&=&\langle h, g \rangle_{A}\\
a_1\ell(h)a_2&=&\ell(a_1ha_2)
\end{eqnarray*}
for any $h,g\in {\mathcal H}$ and any $a_1,a_2\in A$.
\item
The sums
$$
\ell(h)+\ell(h)^*
$$
are $A$-{\it valued semicircular elements}.
\item
In a similar way one can define the full Fock space over any Hilbert $A$-bimodule ${\mathcal H}$, where $A$ is a unital $C^*$-algebra, as well as the associated left creation operators. 
Pimsner showed in [11] that the $C^*$-algebra ${\mathcal T}_{\mathcal{H}}$ generated by $A$ and the 
family of elements $\{\ell(h): h\in{\mathcal H}\}$ subject to the above relations has the universal property. This $C^*$ algebra is called the {\it Toeplitz algebra} over ${\mathcal H}$. 
If $\mathcal{H}$ is {\it full}, which means that $A$ is generated by the set of inner products $\{\langle h, g\rangle_{A}: h,g\in \mathcal{H}\}$, the 
copy of $A$ acting on the left of ${\mathcal F}({\mathcal H})$ is contained in 
$C^*(\ell(h), h\in \mathcal{H})$ and thus
$$
\mathcal{T}_{\mathcal{H}}=C^*(\ell(h), h\in \mathcal{H}).
$$
We will be interested in Toeplitz algebras ${\mathcal T}_{\mathcal H}$ associated with 
Hilbert $A$-bimodules of the form described in (1)-(3), and in the associated 
quotient $C^*$-algebras ${\mathcal O}_{{\mathcal H}}$ which can be identified with 
the {\it Cuntz-Krieger algebras} [2].
\end{enumerate}

We can show that the $C^*$-algebra ${\mathcal T}_{\mathcal{M}}$ generated by the family of standard 
matricially free creation operators is isomorphic to the Toeplitz algebra
${\mathcal T}_{\mathcal{H}}$ specified above. 

\begin{Theorem} 
If ${\mathcal T}_{\mathcal{M}}:=C^*(\wp_{p,q}(u):(p,q)\in \mathpzc{J}, u\in \mathpzc{U})$, 
where each $\wp_{p,q}(u)$ is standard, then 
$$
{\mathcal T}_{\mathcal{M}}\cong {\mathcal T}_{{\mathcal H}},
$$
where 
${\mathcal T}_{\mathcal{H}}$ is the Toeplitz algebra associated with the Hilbert $A$-bimodule ${\mathcal H}$, 
described in (1)-(3), where each ${\mathcal H}_{p,q}(u)$ is a one-dimensional complex Hilbert space.
\end{Theorem}
{\it Proof.}
Using standard matricially free creation operators, we can define general creation operators 
on ${\mathcal M}$ which are associated with elements of ${\mathcal H}$. Namely, let
\begin{eqnarray*}
\wp(h)&=&\sum_{u\in \mathpzc{U}}\sum_{(p,q)\in \mathpzc{J}}h_{p,q}^{c}(u)\wp_{p,q}(u)
\end{eqnarray*}
be associated with the vector $h\in \mathcal{H}$ such that
$h_{p,q}(u)=h_{p,q}^{c}(u)e_{p,q}(u)$ for any $p,q,u$.
Now, let us distinguish subspaces of ${\mathcal M}$ of the form
$$
{\mathcal M}(p)={\mathbb C}\Omega_{p}\oplus \bigoplus_{m=1}^{\infty}
\bigoplus_{\stackrel{p,p_1,\ldots , p_{m}}
{\scriptscriptstyle u_1, \ldots , u_m}
}
{\mathcal H}_{p,p_1}(u_{1})\otimes {\mathcal H}_{p_1,p_2}(u_2)
\otimes \ldots \otimes 
{\mathcal H}_{p_{m-1},p_{m}}(u_{m}),
$$
where $p\in [r]$, and denote by $P_p$ the associated canonical orthogonal projection.
Clearly, these projections commute with each other and $P_1+\ldots +P_r=1$. 
It is easy to see that ${\mathcal T}_{\mathcal{M}}$ coincides with the 
$C^*$-algebra generated by $P_1, \ldots , P_r$ and by the family of operators 
$\{\wp(h), h\in \mathcal{H}\}$, since 
$\wp_{p,q}(u)=P_{p}\wp(h)P_q$. Moreover, it is easy to show that 
the generators $\wp(h)$ and $P_q$ satisfy the following relations:
\begin{eqnarray*}
\alpha\wp(h)+\beta\wp(g)&=&\wp(\alpha h + \beta g)\\
P_{p}\wp(h)P_q&=&\wp(F_{p}hF_q)\\
\wp(h)^*\wp(g)&=&\sigma(\langle h, g\rangle_{A})
\end{eqnarray*}
for any $p,q\in [r],h,g \in {\mathcal H}$ and $\alpha, \beta\in {\mathbb C}$, where
$\sigma: A\rightarrow {\mathcal T}_{{\mathcal M}}$ is the *-homomorphism defined by 
$\sigma(F_q)=P_q$ for any $q$. For instance, 
\begin{eqnarray*}
\wp(h)^*\wp(g)&=&(\sum_{u\in \mathpzc{U}}\sum_{p,q}\overline{h_{p,q}^c(u)}\wp_{p,q}(u)^*)
(\sum_{v\in \mathpzc{U}}\sum_{s,t}g_{s,t}^c(v)\wp_{s,t}(v))\\
&=&
\sum_{u\in \mathpzc{U}}\sum_{p,q}\overline{h_{p,q}^{c}(u)}g_{p,q}^{c}(u)P_q\\
&=&
\sigma(\sum_{u\in \mathpzc{U}}\sum_{p,q}\langle h_{p,q}(u), g_{p,q}(u)\rangle F_{q})\\
&=&
\sigma(\langle h, g \rangle_{A}).
\end{eqnarray*}
The remaining relations are also straightforward. In view of [11, Theorem 3.4], 
the *-homomorphism $\sigma$ has a unique extension 
$\widetilde{\sigma}:{\mathcal T}_{\mathcal{H}}\rightarrow {\mathcal T}_{\mathcal {M}}$ 
that maps $\ell(h)$ to $\wp(h)$, which implies that 
${\mathcal T}_{\mathcal{H}}\cong {\mathcal T}_{\mathcal {M}}$.
\hfill $\blacksquare$\\

\begin{Remark}
{\rm Let us make some remarks referring to Theorem 2.1, its proof and some of its 
consequences.
\begin{enumerate}
\item
In particular, if $h_{p,q}^{c}(u)=\sqrt{b_{p,q}(u)}$, then $\wp(h)$
becomes the sum of matricially free creation operators with covariances 
$b_{p,q}(u)$ of Definition 2.1 and its counterpart $\ell(h)$ in ${\mathcal F}$ is
the canonical creation operator corresponding to the matrix $B(u)$.
In particular, if each $b_{p,q}(u)=d_{p}\geq 0$ and 
$d_1+\ldots +d_r=1$ for any $(p,q)\in \mathpzc{J}=[r]\times [r]$ and 
$u\in \mathpzc{U}$, then $\wp(h)$ is an isometry.
\item
The Hilbert $A$-bimodule ${\mathcal H}$ studied here is finitely generated 
and $A$ is a commutative finite dimensional $C^*$-algebra. In this case, 
let $a_{p,q}$ be the dimension of 
$F_p{\mathcal H}F_q=\bigoplus_{u\in \mathpzc{U}}{\mathcal H}_{p,q}(u)$.
The family $\{\widetilde{e}_{p,q}(u):(p,q)\in \mathpzc{J}, u\in \mathpzc{U}\}$, 
where $\widetilde{e}_{p,q}(u)$ is the canonical injection of $e_{p,q}(u)$ into ${\mathcal H}$,
is the basis of the vector space ${\mathcal H}$ which satisfies conditions
\begin{eqnarray*}
F_{p'}\widetilde{e}_{p,q}(u)F_{q'}&=&\delta_{p',p}\delta_{q,q'}\widetilde{e}_{p,q}(u)\\
\langle \widetilde{e}_{p,q}(u), \widetilde{e}_{p',q'}(u')\rangle_{A}
&=&\delta_{p,p'}\delta_{q,q'}\delta_{u,u'}F_{q}
\end{eqnarray*}
as described by Pimsner [11]. 
\item
With each Hilbert $A$-bimodule ${\mathcal H}$ one can associate the Cuntz-Pimsner algebra 
${\mathcal O}_{{\mathcal H}}$ which is a quotient of ${\mathcal T}_{\mathcal H}$.
Suppose, for simplicity, that we have one array and thus $u$ can be omitted,
and consider sums of the form
$$
\widetilde{e}_{p}=\sum_{q}a_{p,q}\widetilde{e}_{p,q},
$$
where each $a_{p,q}\in \{0,1\}$ and it is assumed that 
the matrix $(a_{p,q})$ has no zero rows and no zero columns. 
The corresponding partial isometries 
$S_{p}=S_{\widetilde{e}_{p}}$ generating the Cuntz-Pimsner algebra 
${\mathcal O}_{{\mathcal H}}$ are classes associated with the creation operators 
$\ell(\widetilde{e}_{p})$. They have orthogonal ranges and satisfy the crucial relation between source projections 
and range projections
\begin{eqnarray*}
S_{p}^*S_{p}&=&\sum_{q}a_{p,q}S_{q}S_{q}^*,
\end{eqnarray*}
thus they generate the Cuntz-Krieger algebra associated with the matrix $(a_{p,q})$.
\item
Using the Hilbert space ${\mathcal M}$ rather than the Hilbert $A$-bimodule ${\mathcal F}$, 
we can associate $\wp_{p}:=\wp(\widetilde{e}_{p})$ 
with each vector $\widetilde{e}_{p}$ and it takes the
form
$$
\wp_p=\sum_{q}a_{p,q}\wp_{p,q}
$$
for any $p$. Each $\wp_{p}$ is a partial isometry with 
a non-zero source projection and non-zero orthogonal range projections. Moreover,
these projections satisfy
$$
\wp_{p}^*\wp_{p}=\sum_{q}a_{p,q}(\wp_{q}\wp_{q}^*+P_{\Omega_{q}}),
$$
where $P_{\Omega_{q}}$ is the orthogonal projection onto ${\mathbb C}\Omega_q$. 
\item
If ${\mathcal O}_{\mathcal{M}}={\mathcal T}_{\mathcal{M}}/{\mathcal I}_{\mathcal{M}}$, 
where ${\mathcal I}_{M}$ is the ideal generated by the set
$$
\{\wp_{p}^*\wp_{p}-\sum_{q}a_{p,q}\wp_{q}\wp_{q}^*: p\in [r], u\in \mathpzc{U}\},
$$
then the corresponding classes $t_{q}=\wp_{p}+\mathcal{I}_{\mathcal{M}}$
satisfy the relation
\begin{eqnarray*}
t_{p}^*t_{p}&=&\sum_{q}a_{p,q}t_{q}t_{q}^*
\end{eqnarray*}
and thus ${\mathcal O}_{\mathcal {M}}\cong {\mathcal O}_{\mathcal {H}}$,
where ${\mathcal O}_{\mathcal{H}}$ is the Cuntz-Krieger 
algebra associated with the matrix $(a_{p,q})$. 
\end{enumerate}}
\end{Remark}

Finally, let us establish a relation between the full Fock space 
${\mathcal F}(\mathcal{H})$ over the Hilbert $A$-bimodule ${\mathcal H}$ 
and the matricially free Fock space ${\mathcal M}$. Clearly, the first object is
a Hilbert space and the second one is a Hilbert $A$-bimodule, but 
both objects can be decomposed into direct sums whose summands 
can be identified. The information about inner products in both structures
is contained in the result stated below. 

\begin{Proposition}
The full Fock space ${\mathcal F}({\mathcal H})$ over the Hilbert $A$-bimodule ${\mathcal H}$, 
equipped with the scalar-valued inner product 
$$
\langle x, y \rangle_{{\mathcal F}} :={\rm Tr}(\langle x, y \rangle_{A})
$$
is unitarily isomorphic to the matricially free Fock space $\mathcal{M}$.
\end{Proposition}
{\it Proof.}
Clearly, $\langle ., .\rangle$ defines a scalar-valued inner product on ${\mathcal F}({\mathcal H})$, with respect to which this space is complete since all norms on $A$ are equivalent.
Moreover, the bijection
\begin{eqnarray*}
\tau(\Omega_{q})&=&F_q \\
\tau(e_{p_1,q_1}(u_1)\otimes \ldots \otimes e_{p_m,q_m}(u_m))&=&\widetilde{e}_{p_1,q_1}(u_1)\otimes_{A}\ldots \otimes_{A}\widetilde{e}_{p_m,q_m}(u_m)
\end{eqnarray*}
extends to a unitary isomorphim of complex Hilbert spaces 
$\tau:{\mathcal M}\rightarrow {\mathcal F}({\mathcal H})$ 
since 
\begin{eqnarray*}
\langle \tau(x), \tau(y) \rangle_{{\mathcal F}}&=&
\sum_{q}{\rm Tr}(\langle \tau(x_{q}),\tau(y_q)\rangle_{A})\\ 
&=&\sum_{q}{\rm Tr}(\langle x_q,y_q\rangle F_q)\\
&=&\sum_{q}\langle x_q, y_q\rangle \\
&=&
\langle x , y \rangle
\end{eqnarray*}
for any $x=(x_q),y=(y_q)\in {\mathcal M}=\bigoplus_{q}{\mathcal M}_{q}$, where 
we use the fact that each subspace ${\mathcal M}_{q}$ is mapped by $\tau$ onto 
$$
{\mathcal F}({\mathcal H})F_q={\mathbb C}F_q\oplus \bigoplus_{m=1}^{\infty} 
\bigoplus_{q_1, \ldots, q_{m}}F_{q_1}{\mathcal H}F_{q_2}\otimes_{A}
F_{q_2}{\mathcal H}F_{q_3}\otimes_{A}
\ldots \otimes_{A} F_{q_{m}}{\mathcal H}F_{q},
$$
and that
${\mathcal F}({\mathcal H})=\bigoplus_{q}{\mathcal F}({\mathcal H})F_q$.
Thus, our proof is completed.
\hfill $\blacksquare$\\

In view of the relation between ${\mathcal M}$ and ${\mathcal F}$, we can 
say that ${\mathcal M}$ is the Hilbert space which is naturally 
associated with ${\mathcal F}$, treated as a Hilbert $A$-bimodule.

\section{Matricial circular systems}

In this section, we will introduce matricial circular systems of 
operators which arise naturally as matricial analogs of circular operators. 
The framework of Hilbert $A$-bimodules discussed above makes them also very 
natural and quickly leads to their realizations as operator-valued matrices.

\begin{Definition}
{\rm By {\it matricial circular operators} we shall understand operators on 
${\mathcal M}$ of the form
$$
\zeta_{p,q}(u)=\wp_{p,q}(u')+\wp^{*}_{q,p}(u''),
$$
where $u',u''$ are copies of $u\in \mathpzc{U}$ and $p,q\in [r]$
and $B(u)=B(u')=B(u'')$.
The corresponding families of arrays of operators 
will be called {\it matricial circular systems}.
}
\end{Definition}

In our study of distributions of matricial systems of operators in our previous works, 
we considered the family $\{\Psi_1, \ldots , \Psi_r\}$ of vector states on $B({\mathcal M})$ 
of the form
$$
\Psi_q(a)=\langle a\Omega_q, \Omega_q\rangle,
$$
from which we constructed the array $(\Psi_{p,q})$ by setting $\Psi_{p,q}=\Psi_q$
for any $p,q$. This array is needed when considering the concepts of matricial freeness and
symmetric matricial freeness. In order to express matricial circular systems 
as operator-valued matrices, we need to recall a few notions. 
\begin{enumerate} 
\item
When speaking of operator-valued matrices, we deal with 
a unital *-algebra ${\mathcal A}$ and a state $\varphi$ on ${\mathcal A}$.
We then obtain a *-probability space $({\mathcal A}, \varphi)$. 
The *-algebra of interest is then the algebra of matrices 
$M_{r}({\mathcal A})\cong {\mathcal A}\otimes M_{r}({\mathbb C})$ with the natural involution 
$$
(a\otimes e(p,q))^{*}=a^{*}\otimes e(q,p)
$$
for any $a\in {\mathcal A}$ and $p,q\in [r]$, where $\{e(p,q): p,q\in [r]\}$
is the system of matrix units. We denote by
$\Phi_1, \ldots , \Phi_r$ the states on $M_{r}({\mathcal A})$ of the form
$$
\Phi_k=\varphi\otimes \psi_k
$$
for any $k\in [r]$, where $\psi_k(b)=\langle b e(k), e(k)\rangle$ and
$\{e(1), \ldots, e(r)\}$ is the canonical orthonormal basis in ${\mathbb C}^{r}$. 
If ${\mathcal A}$ is a $C^{*}$-algebra, we can work 
in the category of $C^{*}$-probability spaces. 
\item
Let us suppose that in the given $C^{*}$-probability space 
$({\mathcal A}, \varphi)$ we have a family of free creation operators,
$
\{\ell(p,q,u): p,q\in [r], \;u\in \mathpzc{U}\},
$
which is *-free with respect to $\varphi$, for which
$$
\ell(p,q,u)^*\ell(p',q',u')=\delta_{p,p'}\delta_{q,q'}\delta_{u,u'}b_{p,q}(u),
$$ 
for any $p,q,u,p',q',u'$, where the corresponding matrices $B(u)=(b_{p,q}(u))$ 
consist of non-negative numbers. The number $b_{p,q}(u)$ will be called the 
{\it covariance} of $\ell(p,q,u)$. A family of arrays of the above form will be 
called a system of free creation operators. If it contains only standard free 
creation operators, we will say that this system is standard.
\item
By a {\it generalized circular element} of a *-probability space $({\mathcal A}, \varphi)$ 
we understand an element whose *-distribution agrees with the *-distribution of 
the sum of the form 
$$
c=\ell_1 + \ell_2^{*},
$$ 
where $\ell_1$ and $\ell_2$ are free creation operators
with covariances $\alpha>0$ and $\beta>0$, respectively, which are *-free with respect to $\varphi$.
If we need to be more specific, we will call the above element an $(\alpha, \beta)$-{\it circular element}. 
In particular, when $\alpha=\beta$, it is called the {\it circular element} with covariance $\alpha$. If $\alpha=1$, the circular element will be called {\it standard}.
\item
By a {\it (generalized) circular system} in $({\mathcal A}, \varphi)$ 
we shall understand the family 
$$
\{g(p,q,u): p,q \in[r], u\in \mathpzc{U}\}
$$
where each $g(p,q,u)=\ell(p,q,u')+\ell^*(q,p,u'')$ (for $u'\neq u''$) 
is a (generalized) circular element and the whole family is *-free with respect to $\varphi$. 
\end{enumerate}

\begin{Proposition}
Let $\{\ell(p,q,u):(p,q)\in \mathpzc{J}, u\in \mathpzc{U}\}$ 
be the family of *-free creation operators with respect to the state $\varphi$
with covariances given by matrices $B(u)$. 
Then the joint *-distribution of the operator-valued matrices
$$
\ell(p,q,u)\otimes e(p,q)
$$
in the state $\Phi_k$, where $p,q,k\in [r]$ and $u\in \mathpzc{U}$, 
agrees with the joint *-distribution of the corresponding matricially
free creation operators $\wp_{p,q}(u)$ in the state $\Psi_k$, respectively.
\end{Proposition}
{\it Proof.}
Denoting $\ell_{p,q}(u)=\ell(\widetilde{e}_{p,q}(u))$ and using the relations
between the generators of $\mathcal{T}_{{\mathcal H}}$, we obtain
$$
\ell_{p,q}(u)^*\ell_{p',q'}(u')=\delta_{p,p'}\delta_{q,q'}\delta_{u,u'}b_{p,q}(u)F_{q},
$$
$$
F_{p'}\ell_{p,q}(u)F_{q'}=\delta_{q,q'}\delta_{p,p'}\ell_{p,q}(u).
$$
Analogous relations hold if we replace $\ell_{p,q}(u)$ by $\ell(p,q,u)\otimes e(p,q)$,
where $\ell(p,q,u)$ has covariance $b_{p,q}(u)$, and we write $F_{q}$ as $1\otimes e(q,q)$.
This means that the joint *-distribution of the  
operator-valued matrices $\ell(p,q,u)\otimes e(p,q)$ in the states 
$\Phi_k$ agree with the joint *-distributions 
of the operators $\ell_{p,q}(u)$ in the states 
$\widetilde{\Phi}_k(.):=\langle F_{k}, .\,F_{k}\rangle_{\mathcal{F}}$, 
respectively. In turn, by the proofs of Theorem 2.1 and Proposition 2.1,
the latter agree with the joint *-distributions 
of the operators $\wp_{p,q}(u)$ in the states $\Psi_k$, which completes the proof.
\hfill $\blacksquare$\\

One can study these *-distributions, using the language of conditional expectations 
defined by a positive matrix $(\eta_{u,w})$ of linear maps on $A$.
\begin{enumerate}
\item
Using the relations between generators of the Toeplitz algebra ${\mathcal T}_{{\mathcal H}}$, 
we get 
\begin{eqnarray*}
\ell(u)^*F_p\ell(w)&=&\ell(u)^*\ell(F_pw)\\
&=&\langle h(u), F_ph(w)\rangle_{A}\\
&=&
\sum_{q}\langle h_{p,q}(u), h_{p,q}(w)\rangle F_q\\
&=&
\delta_{u,w}\sum_{q}b_{p,q}(u)F_q\\
\end{eqnarray*}
where $\ell(u)=\ell(h(u))$ and $\ell(w)=\ell(h(w))$ for some $h(u),h(w)\in {\mathcal H}$.
\item
We can write these relations in the form
$$
\ell(u)^*F_p \ell(w)=\eta_{u,w}(F_{p}),
$$
where $\eta_{u,w}:A\rightarrow A$ is a linear map given by 
$$
\eta_{u,w}(F_p)=\delta_{u,w}\sum_{q}b_{p,q}(u)F_q
$$
for any $u,w\in \mathpzc{U}$ and $p\in [r]$. 
\item
We know from [12] that positive matrices of linear maps implemented in this fashion 
define a conditional expectation 
$E: {\mathcal T}_{{\mathcal H}}\rightarrow A$ 
by setting $E(F_p)=F_p$ and $E(F_pxF_q)=F_pE(x)F_q$
for any $p,q\in [r]$ and $x\in {\mathcal T}_{{\mathcal H}}$, and requiring that 
$E(W)=0$ for all words $W$ in $\ell(u)^*, \ell(u), F_p$, where 
$u\in \mathpzc{U}$ and $p\in [r]$, which are not reducible to 
an element of $A$ by the above relations.  
\item
Identifying the generators $\ell(u)$ of the Toeplitz algebra ${\mathcal T}_{{\mathcal H}}$
with the operator-valued matrices 
$$
L(u)=\sum_{p,q}\ell(p,q,u)\otimes e(p,q),
$$
where $u\in \mathpzc{U}$ and each $\ell(p,q,u)\in \mathcal{A}$ has covariance $b_{p,q}(u)$,
we can identify the conditional expectation on the Toeplitz algebra described above 
with the conditional expectation $E:M_{r}({\mathcal A})\rightarrow A$ of the form
$$
E\left(\sum_{p,q}a(p,q)\otimes e(p,q)\right)=\sum_{q}\varphi(a(q,q))\otimes e(q,q)
$$
where $a(p,q)\in \mathcal{A}$ for any $p,q\in [r]$ and $\varphi$ is the vacuum expectation 
on ${\mathcal A}$. Then the *-distributions of Proposition 3.1 are recovered by computing 
$E(F_kxF_k)$ for suitable $x\in M_{r}({\mathcal A})$ and any $k\in [r]$. 
Operators $L(u)$ (in fact, even more general ones) were studied in [13].
\end{enumerate}
  
\begin{Corollary}
With the above notations, suppose that each $g(p,q,u)$ is $(\alpha, \beta)$-circular, where $\alpha=b_{p,q}(u)$and $\beta=b_{q,p}(u)$. Then the joint *-distributions of the operators
$$
g(p,q,u)\otimes e(p,q)
$$
in the states $\Phi_k$, where $p,q,k\in [r]$ and $u\in \mathpzc{U}$, 
agree with the joint *-distributions of the corresponding matricial 
circular operators $\zeta_{p,q}(u)$ in the states $\Psi_k$, respectively.
\end{Corollary}
{\it Proof.}
It suffices to use the expression $g(p,q,u)=\ell(p,q,u')+\ell(q,p,u'')^*$
for any $p,q,u$, and Proposition 3.1 to obtain this assertion.\hfill $\blacksquare$.\\

Note that these operator-valued matrices are the generators of free group factors used by Voiculescu in his proof of free group factors isomorphisms [16, Theorem 3.3].

\section{Asymptotics of random blocks}

Matricial semicircular (symmetrized circular) operators give operatorial realizations of the limit joint 
distributions (*-distributions) under partial traces of symmetric random blocks of Hermitian 
(non-Hermitian) GRM with i.b.i.d. complex entries, respectively [6,7].
Let us prove an analogous result for (non-symmetric) blocks of GRM and 
matricial circular operators.

Let us recall basic definitions needed in our study of the asymptotics of blocks.
\begin{enumerate}
\item
Consider the partition of $[n]:=\{1,2, \ldots, n\}$ 
into disjoint non-empty intervals 
$$
[n]=N_{1}\cup \ldots \cup N_{r}
$$
whose dependence on $n$ is supressed, and let $n_q=|N_q|$. Further, assume that
$$
\lim_{n\rightarrow \infty}\frac{n_{q}}{n}= d_{q}\geq 0\;\;{\rm for}\;{\rm any}\; q,
$$
and let $D=(d_1, d_2, \ldots, d_r)$ be the so-called {\it dimension matrix}.
\item
By {\it partial traces} we understand normalized traces composed with classical expectations of the form 
$$
\tau_q(n)=\frac{1}{n_q}\left(E\otimes {\rm tr}_{q}(n)\right)
$$
where $q\in [r]$ and ${\rm tr}_{q}(n)$ is the trace over the set
of basis vectors of ${\mathbb C}^{n}$ indexed by $N_{q}$ for any $q$. 
\item
We would like to study the asymptotics of blocks of i.b.i.d. $n\times n$ Gaussian random matrices
$Y(u,n)=(Y_{i,j}(u,n))$ of the form
$$
S_{p,q}(u,n)=\sum_{i,j\in N_p\times N_q}Y_{i,j}(u,n)\otimes e_{i,j}(n),
$$
where $p,q\in [r]$ and $\{e_{i,j}(n):1\leq i,j \leq n\}$ is the system of matrix units in $M_{n}({\mathbb C})$, 
under partial traces.
We clearly have 
$$
S_{p,q}(u,n)=D_{p}Y(u,n)D_{q},
$$
where $D_{q}$ is the $n\times n$ diagonal matrix such that $(D_{q})_{k,k}=1$ whenever $k\in N_{q}$ and 
its remaining entries are zeros.
\item
The symmetric blocks studied in [7] are constructed from these blocks in the natural way, namely
$$
T_{p,q}(u,n)=
\left\{
\begin{array}{ll}
S_{p,q}(u,n)+S_{q,p}(u,n)&{\rm if} \;\;p< q\\
S_{q,q}(u,n) & {\rm if} \;\;p=q
\end{array}
\right.
$$
for any $p,q,u,n$. Note that we identify blocks and symmetric blocks with suitable subblocks of 
$n\times n$ matrices. 
\end{enumerate}

\begin{Theorem}
Let $\{Y(u,n): u\in \mathpzc{U}\}$ be a family of GRM with i.b.i.d. 
complex entries, each with blocks $(S_{p,q}(u,n))$ and symmetric block covariance matrix $V(u)$.
Then
$$
\lim_{n\rightarrow \infty}S_{p,q}(u,n)=\zeta_{p,q}(u)
$$
for any $1\leq p,q\leq r$ and $u\in \mathpzc{U}$, where each array 
$(\zeta_{p,q}(u))$ is associated with matrix $B(u)=DV(u)$ and
convergence is in mixed *-moments under partial traces.
\end{Theorem}
{\it Proof.}
In the case of non-degenerate dimension matrix $D$, 
one can use [12, Theorem 4.1] which implies that
$$
\lim_{n\rightarrow \infty}D_{p}Y(u,n)D_{q}=F_pg(u)F_q
$$
in the sense of mixed *-moments under $E$, where 
$g(u)=\ell(u')+\ell(u'')^*$, and then 
reproduce the limit mixed *-moments under $\tau_{q}(n)$ 
by taking the compressions of those on the RHS by $F_q$.
A proof using scalar-valued free probability, similar to that in [7, Theorem 6.1],
can also be given. Namely, in the i.i.d. case and non-degenerate $D$ 
we use asymptotic *-freeness of the family $\{Y(u,n):u\in \mathpzc{U}\}$ and its asymptotic *-freeness
from $\{D_1, \ldots , D_r\}$ under $\tau(n)$, which gives
$$
\lim_{n\rightarrow \infty}D_pY(u,n)D_q=P_p\zeta(u)P_q
$$
in the sense of convergence of mixed *-moments under $\tau(n)$ to mixed *-moments under $\Psi$, 
where $\zeta(u)=\wp(u')+\wp(u'')^*$. The i.b.i.d. case is then obtained by rescaling blocks. 
If $D$ is singular, we can apply a limiting procedure, as in [7, Theorem 6.1]. The details are left to the reader.
\hfill $\blacksquare$\\

Of course, this immediately gives the analogous result for symmetric blocks 
(that one already appeared in [7]).

\begin{Corollary}
Under the assumptions of Theorem 4.1, it holds that
$$
\lim_{n\rightarrow \infty}T_{p,q}(u,n)=\eta_{p,q}(u)
$$
for any $1\leq p,q\leq r$ and $u\in \mathpzc{U}$, where convergence is in the sense 
of mixed *-moments under partial traces.
\end{Corollary}

Knowing that the limit mixed *-distributions of blocks 
under partial traces are given by *-distributions of matricial circular operators, we 
are going to study the latter in more detail.

\section{Joint *-distributions}

The combinatorics of mixed *-moments of matricial circular operators under states $\Phi_q$ is based on non-crossing colored pair partitions and is similar to that for symmetrized matricial circular operators [7]. Its main feature is that it suffices to color the blocks, including the imaginary block colored by $q$, and assign covariances to these blocks in a matricial way. This combinatorics can be viewed as a special case of the more general combinatorics of mixed *-moments of operator-valued circular elements, which can 
be reduced to that for operator-valued semicircular elements [14]. 

If $\pi$ is a noncrossing pair-partition of the set $[m]$, where $m$ is an even positive integer, 
which is denoted $\pi\in \mathcal{NC}_{m}^{2}$, the set 
$$
{\mathcal B}(\pi)=\{V_1, V_2, \ldots , V_s\}
$$ 
is the set of its blocks, where $m=2s$. If $V_{i}=\{l(i),r(i)\}$ and $V_{j}=\{l(j),r(j)\}$ are two blocks of
$\pi$ with left legs $l(i)$ and $l(j)$ and right legs $r(i)$ and $r(j)$, respectively, then
$V_i$ is {\it inner} with respect to $V_j$ if $l(j)<l(i)<r(i)<r(j)$.
In that case $V_j$ is {\it outer} with respect to $V_i$.
It is the {\it nearest outer block} of $V_i$ if there is no block $V_k=\{l(k),r(k)\}$
such that $l(j)<l(k)<l(i)<r(i)<r(k)<r(j)$.
It is easy to see that the nearest outer block, if it exists, 
is unique, and we write in this case $V_j=o(V_i)$. If $V_i$ does not have an outer block, 
we set $o(V_i)=V_0$, where $V_{0}=\{0,m+1\}$ is the additional block called {\it imaginary}. The partition of the set $\{0,1, \ldots , m+1\}$ consisting of the blocks of $\pi$ and of the imaginary block will be denoted by $\widehat{\pi}$.

\begin{Definition}
{\rm 
The tuple $(e(p_1,q_1), \ldots , e(p_m,q_m))$ of matrix units in $M_{r}({\mathbb C})$ 
will be called {\it cyclic} if 
$$
q_1=p_2, q_2=p_3, \ldots , q_{m}=p_1.
$$
We will say that $\pi\in \mathcal{NC}_{m}$ is {\it adapted to} 
$(e(p_1,q_1), \ldots , e(p_m,q_m))$ if this tuple and all tuples associated with 
the blocks of $\pi$ are cyclic. The set of all noncrossing pair partitions adapted
to the tuple $(e(p_1,q_1), \ldots , e(p_m,q_m))$, where $p_1=q_m$, will be denoted by 
$\mathcal{NC}_{m}(e(p_1,q_1), \ldots ,e(p_m,q_m))$.}
\end{Definition}

\begin{Example}
{\rm For a tuple of off-diagonal matrix units $(e(p,q),e(q,p), e(p,q), e(q,p))$,
there are three noncrossing partitions which are 
adapted to it, namely
$$
\{\{1,2,3,4\}\},\;\;\{\{1,4\}, \{2,3\}\}\;\; {\rm and}\;\;\{\{1,2\}, \{3,4\}\}
$$ 
since only the tuples of matrix units associated to these partitions are cyclic.
}
\end{Example}

\begin{Definition}
{\rm 
If we are given $a_j=c_j\otimes e(p_j,q_j)\in M_{r}({\mathcal A})$ 
for $j=1, \ldots, m$, where $c_j\in {\mathcal A}$, we will denote by  
$$
\mathcal{NC}_{m}(a_1, \ldots, a_m)
$$ 
the set of noncrossing partitions of $[m]$ which are adapted to $(e(p_1,q_1), \ldots , e(p_m,q_m))$.
These partitions will be called {\it adapted to the tuple} $(a_1, \ldots , a_m)$.
Its subset consisting of pair partitions will be denoted 
$\mathcal{NC}^{2}_{m}(a_1, \ldots, a_m)$. }
\end{Definition}

The mixed *-moments of the circular operators $\zeta_{p,q}(u)$ under the states $\Phi_k$
can be expressed in terms of noncrossing pair partitions which are adapted to the associated tuples.
The contribution from each such partition $\pi$ is a product of covariances. In order to specify this 
contribution, we shall color the blocks of $\pi$ by numbers from the set $[r]$.
If $\pi\in \mathcal{NC}_{m}^{2}$, where $m$ is even, we denote by $\mathpzc{F}_{r}(\pi)$ the set of all mappings 
$f:{\mathcal B}(\pi)\rightarrow [r]$ called {\it colorings}. 
By a {\it colored noncrossing pair partition} we then understand a pair $(\pi,f)$, where 
$\pi \in \mathcal{NC}_{m}^{2}$ and $f\in \mathpzc{F}_{r}(\pi)$. The set of pairs 
$$
{\mathcal B}(\pi,f)=\{(V_1,f), (V_2,f), \ldots , (V_k,f)\}
$$ 
will play the role of the set of its blocks. We will always assume that also the imaginary block is colored by a number from the set $[r]$ and thus we can speak of a coloring of $\widehat{\pi}$.  
Examples of colored noncrossing pair partitions are given in Fig.1.

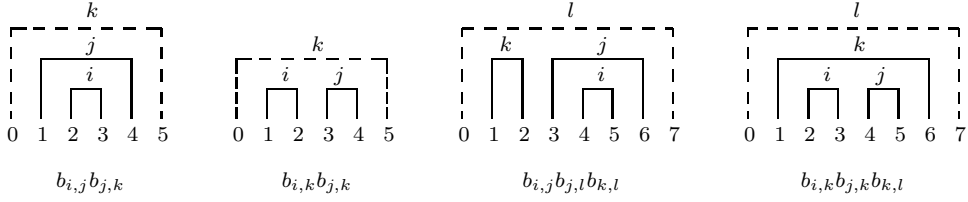
\begin{figure}
\unitlength=1mm
\special{em:linewidth 0.4pt}
\linethickness{0.4pt}
\begin{picture}(120.00,35.00)(30.00,0.00)
\put(31.00,10.00){\line(0,1){8.00}}
\put(35.00,10.00){\line(0,1){4.00}}
\put(39.00,10.00){\line(0,1){4.00}}
\put(43.00,10.00){\line(0,1){8.00}}

\put(27.00,10.00){\line(0,1){1.00}}
\put(27.00,12.00){\line(0,1){1.50}}
\put(27.00,15.00){\line(0,1){1.50}}
\put(27.00,18.00){\line(0,1){1.50}}
\put(27.00,21.00){\line(0,1){1.00}}

\put(47.00,10.00){\line(0,1){1.00}}
\put(47.00,12.00){\line(0,1){1.50}}
\put(47.00,15.00){\line(0,1){1.50}}
\put(47.00,18.00){\line(0,1){1.50}}
\put(47.00,21.00){\line(0,1){1.00}}

\put(26.50,7.00){$\scriptstyle{0}$}
\put(30.50,7.00){$\scriptstyle{1}$}
\put(34.50,7.00){$\scriptstyle{2}$}
\put(38.50,7.00){$\scriptstyle{3}$}
\put(42.50,7.00){$\scriptstyle{4}$}
\put(46.50,7.00){$\scriptstyle{5}$}

\put(37.00,15.00){$\scriptstyle{i}$}
\put(37.00,19.00){$\scriptstyle{j}$}
\put(37.00,23.50){$\scriptstyle{k}$}

\put(27.00,22.00){\line(1,0){2.00}}
\put(30.00,22.00){\line(1,0){1.50}}
\put(33.00,22.00){\line(1,0){1.50}}
\put(36.00,22.00){\line(1,0){1.50}}
\put(39.00,22.00){\line(1,0){1.50}}
\put(42.00,22.00){\line(1,0){1.50}}
\put(45.00,22.00){\line(1,0){2.00}}

\put(31.00,18.00){\line(1,0){12.00}}
\put(35.00,14.00){\line(1,0){4.00}}
\put(61.00,10.00){\line(0,1){4.00}}
\put(65.00,10.00){\line(0,1){4.00}}
\put(69.00,10.00){\line(0,1){4.00}}
\put(73.00,10.00){\line(0,1){4.00}}

\put(57.00,10.00){\line(0,1){1.50}}
\put(57.00,12.00){\line(0,1){1.50}}
\put(57.00,14.00){\line(0,1){1.50}}
\put(57.00,16.00){\line(0,1){2.00}}

\put(77.00,10.00){\line(0,1){1.50}}
\put(77.00,12.00){\line(0,1){1.50}}
\put(77.00,14.00){\line(0,1){1.50}}
\put(77.00,16.00){\line(0,1){2.00}}

\put(56.50,7.00){$\scriptstyle{0}$}
\put(60.50,7.00){$\scriptstyle{1}$}
\put(64.50,7.00){$\scriptstyle{2}$}
\put(68.50,7.00){$\scriptstyle{3}$}
\put(72.50,7.00){$\scriptstyle{4}$}
\put(76.50,7.00){$\scriptstyle{5}$}

\put(63.00,15.00){$\scriptstyle{i}$}
\put(70.00,15.00){$\scriptstyle{j}$}
\put(67.00,19.00){$\scriptstyle{k}$}

\put(57.00,18.00){\line(1,0){2.00}}
\put(60.00,18.00){\line(1,0){1.50}}
\put(63.00,18.00){\line(1,0){1.50}}
\put(66.00,18.00){\line(1,0){1.50}}
\put(69.00,18.00){\line(1,0){1.50}}
\put(72.00,18.00){\line(1,0){1.50}}
\put(75.00,18.00){\line(1,0){2.00}}

\put(61.00,14.00){\line(1,0){4.00}}
\put(69.00,14.00){\line(1,0){4.00}}
\put(91.00,10.00){\line(0,1){8.00}}
\put(95.00,10.00){\line(0,1){8.00}}
\put(99.00,10.00){\line(0,1){8.00}}
\put(103.00,10.00){\line(0,1){4.00}}
\put(107.00,10.00){\line(0,1){4.00}}
\put(111.00,10.00){\line(0,1){8.00}}

\put(87.00,10.00){\line(0,1){1.00}}
\put(87.00,12.00){\line(0,1){1.50}}
\put(87.00,15.00){\line(0,1){1.50}}
\put(87.00,18.00){\line(0,1){1.50}}
\put(87.00,21.00){\line(0,1){1.00}}

\put(115.00,10.00){\line(0,1){1.00}}
\put(115.00,12.00){\line(0,1){1.50}}
\put(115.00,15.00){\line(0,1){1.50}}
\put(115.00,18.00){\line(0,1){1.50}}
\put(115.00,21.00){\line(0,1){1.00}}

\put(86.50,7.00){$\scriptstyle{0}$}
\put(90.50,7.00){$\scriptstyle{1}$}
\put(94.50,7.00){$\scriptstyle{2}$}
\put(98.50,7.00){$\scriptstyle{3}$}
\put(102.50,7.00){$\scriptstyle{4}$}
\put(106.50,7.00){$\scriptstyle{5}$}
\put(110.50,7.00){$\scriptstyle{6}$}
\put(114.50,7.00){$\scriptstyle{7}$}

\put(92.00,19.00){$\scriptstyle{k}$}
\put(105.00,15.00){$\scriptstyle{i}$}
\put(105.00,19.00){$\scriptstyle{j}$}
\put(101.00,23.50){$\scriptstyle{l}$}

\put(87.00,22.00){\line(1,0){1.50}}
\put(90.00,22.00){\line(1,0){1.50}}
\put(93.00,22.00){\line(1,0){1.50}}
\put(96.00,22.00){\line(1,0){1.50}}
\put(99.00,22.00){\line(1,0){1.50}}
\put(102.00,22.00){\line(1,0){1.50}}
\put(105.00,22.00){\line(1,0){1.50}}
\put(108.00,22.00){\line(1,0){1.50}}
\put(111.00,22.00){\line(1,0){1.50}}
\put(114.00,22.00){\line(1,0){1.00}}

\put(99.00,18.00){\line(1,0){12.00}}
\put(103.00,14.00){\line(1,0){4.00}}
\put(91.00,18.00){\line(1,0){4.00}}
\put(129.00,10.00){\line(0,1){8.00}}
\put(133.00,10.00){\line(0,1){4.00}}
\put(137.00,10.00){\line(0,1){4.00}}
\put(141.00,10.00){\line(0,1){4.00}}
\put(145.00,10.00){\line(0,1){4.00}}
\put(149.00,10.00){\line(0,1){8.00}}

\put(125.00,10.00){\line(0,1){1.00}}
\put(125.00,12.00){\line(0,1){1.50}}
\put(125.00,15.00){\line(0,1){1.50}}
\put(125.00,18.00){\line(0,1){1.50}}
\put(125.00,21.00){\line(0,1){1.00}}

\put(153.00,10.00){\line(0,1){1.00}}
\put(153.00,12.00){\line(0,1){1.50}}
\put(153.00,15.00){\line(0,1){1.50}}
\put(153.00,18.00){\line(0,1){1.50}}
\put(153.00,21.00){\line(0,1){1.00}}

\put(124.50,7.00){$\scriptstyle{0}$}
\put(128.50,7.00){$\scriptstyle{1}$}
\put(132.50,7.00){$\scriptstyle{2}$}
\put(136.50,7.00){$\scriptstyle{3}$}
\put(140.50,7.00){$\scriptstyle{4}$}
\put(144.50,7.00){$\scriptstyle{5}$}
\put(148.50,7.00){$\scriptstyle{6}$}
\put(152.50,7.00){$\scriptstyle{7}$}

\put(135.00,15.00){$\scriptstyle{i}$}
\put(142.00,15.00){$\scriptstyle{j}$}
\put(139.00,19.00){$\scriptstyle{k}$}
\put(139.00,23.50){$\scriptstyle{l}$}

\put(125.00,22.00){\line(1,0){1.25}}
\put(127.67,22.00){\line(1,0){1.25}}
\put(130.33,22.00){\line(1,0){1.25}}
\put(132.98,22.00){\line(1,0){1.25}}
\put(135.63,22.00){\line(1,0){1.25}}
\put(138.30,22.00){\line(1,0){1.25}}
\put(140.95,22.00){\line(1,0){1.25}}
\put(143.60,22.00){\line(1,0){1.25}}
\put(146.25,22.00){\line(1,0){1.25}}
\put(149.00,22.00){\line(1,0){1.25}}
\put(151.70,22.00){\line(1,0){1.30}}

\put(129.00,18.00){\line(1,0){20.00}}
\put(133.00,14.00){\line(1,0){4.00}}
\put(141.00,14.00){\line(1,0){4.00}}
\put(33.00,1.00){$\scriptstyle{b_{i,j}b_{j,k}}$}
\put(63.00,1.00){$\scriptstyle{b_{i,k}b_{j,k}}$}
\put(95.00,1.00){$\scriptstyle{b_{i,j}b_{j,l}b_{k,l}}$}
\put(132.00,1.00){$\scriptstyle{b_{i,k}b_{j,k}b_{k,l}}$}

\end{picture}
\caption{Colored noncrossing pair partitions}
\end{figure}

For any given $r\times r$ covariance matrix $B(u)$, there is a natural way to assign its entries to the blocks of noncrossing pair partitions colored by the set $[r]$. By multiplicativity over the blocks, we can
then define the associated functions on the set of noncrossing pair partitions.

\begin{Definition}
{\rm Let a covariance matrix $B(u)=(b_{p,q}(u))\in M_{r}({\mathbb R})$ be given for any $u\in \mathpzc{U}$.
For any $\pi\in \mathcal{NC}_{m}^{2}$ and $f\in \mathpzc{F}_{r}(\pi)$, let
$$
b_{q}(\pi, f)=\prod_{k=1}^{s}b_{q}(V_{k},f)
$$
where 
$$
b_{q}(V_k,f)=b_{s,t}(u),
$$
whenever 
$V_k=\{i,j\}$ is colored by $s$, its nearest outer block $o(V_k)$ is colored by $t$ and $u_i=u_j=u$ and
we assume that the imaginary block is colored by $q\in [r]$, and otherwise we set $b_{q}(V_k,f)=0$. }
\end{Definition}

It remains to determine which colorings are natural for circular operators. It is convenient to introduce the following definition.
\begin{Definition}
{\rm Let $\pi \in \mathcal{NC}_{m}^{2}(a_1, \ldots, a_m)$, where 
$a_j=(c_j(u_j)\otimes e(p_j,q_j))^{\epsilon_j}$ and $\epsilon_j\in \{1,*\}$ for $j\in [m]$ and 
$m$ is even. A coloring $f:{\mathcal B}(\pi)\rightarrow [r]$ will be called {\it adapted} to $(a_1, \ldots , a_m)$ if
$$
f(V_k)=
\left\{
\begin{array}{ll}
p & {\rm if} \;(\epsilon_i,\epsilon_j)=(*,1)\\
q & {\rm if} \;(\epsilon_i,\epsilon_j)=(1,*)
\end{array}
\right.
$$ 
whenever $V_k=\{i,j\}$ is a block and $(p_i,q_i)=(p_j,q_j)=(p,q)$.}
\end{Definition}

\begin{Lemma}
With the above notations, 
let $a_j=\zeta_{p_j,q_j}^{\epsilon_j}(u_j)$, where $p_j,q_j\in [r]$, $u_j\in \mathpzc{U}$ and 
$\epsilon_j\in \{1,*\}$ for $j\in [m]$ and $m\in {\mathbb N}$. Then
$$
\Phi_{q}(a_1\ldots a_m)=\sum_{\pi\in \mathcal{NC}_{m}^{2}
(a_1, \ldots , a_m)}b_{q}(\pi,f)
$$
where $f$ is the unique coloring of $\pi$ which is adapted to 
$(a_1, \ldots, a_m)$.
\end{Lemma}
{\it Proof.}
We adopt the convention that if $m$ is odd, then all corresponding sets of adapted 
pair partitions are empty, in which case both sides of the above formula vanish.
Thus, we only need to compute the mixed *-moments
$$
\Phi_{q}(\zeta_{p_1,q_1}^{\epsilon_1}(u_1)\ldots \zeta_{p_m,q_m}^{\epsilon_m}(u_m))
$$
where $\epsilon_1, \ldots , \epsilon_m\in \{*,1\}$ and $m$ is even.
These mixed *-moments can be identified with operator-valued *-moments of the form
$$
E(F_{q}(F_{p_1}g(u_1)F_{q_1})^{\epsilon_1}\ldots (F_{p_m}g(u_m)F_{q_m})^{\epsilon_m}F_q),
$$
respectively, where 
$$
g(u)=\ell(u')+\ell(u'')^*
$$
are $A$-{\it valued circular elements} with respect to
the conditional expectation $E$ defined in Section 3. In fact,  
the operator-valued moment is equal to the scalar-valued one multiplied by $F_q$.
In fact, it is equal to the sum of mixed *-moments in 
creation and annihilation operators labelled by various $u_k'$ and $u_k''$, interlaced with
diagonal matrices. These moments do not vanish if the tuple 
$((p_1,q_1,\epsilon_1), \ldots, (p_m,q_m,\epsilon_m))$ 
defines a non-crossing pair partition $\pi$, for which
$$
u_i=u_j,\;\;\epsilon_i\neq \epsilon_j\;\;{\rm and}\;\;(p_i,q_i)=(p_j,q_j)
$$
whenever $\{i,j\}$ is a block, which implies that $\pi\in \mathcal{NC}_{m}^{2}(a_1, \ldots ,a_m)$, 
the condition $u_i=u_j$ being included in our definition of $b_{q}(\pi,f)$.
The computation of these moments is based on repeated use
of the relation 
$$
F_q\ell(u)^*F_p \ell(u)F_q=b_{p,q}(u)F_q
$$
for any $p,q,u$, starting from blocks of largest depths.
We need to check how the colors are assigned to blocks of such $\pi$ in order to 
reproduce the appropriate product of covariances. Let $V=\{i,j\}$ be a block and let 
$o(V)=\{o(i),o(j)\}$ be its nearest outer block. This notation also
applies to the blocks whose nearest outer block is the imaginary block to which
we assign the pair $(q,q)$. We have two types of blocks: 
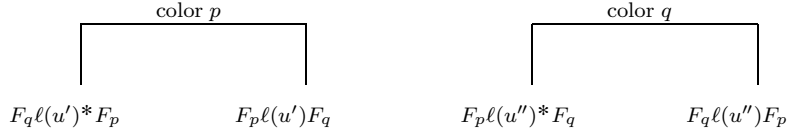
\unitlength=1mm
\special{em:linewidth 0.4pt}
\linethickness{0.4pt}
\begin{picture}(120.00,35.00)(-8.00,-5.00)
\put(20.00,10.00){\line(0,1){8.00}}
\put(50.00,10.00){\line(0,1){8.00}}
\put(20.00,18.00){\line(1,0){30.00}}

\put(30.00,19.00){$\scriptstyle{{\rm color}\;p}$}
\put(10.50,5.00){$\scriptstyle{F_q\ell(u')^*F_p}$}
\put(40.50,5.00){$\scriptstyle{F_p\ell(u')F_q}$}
\put(80.00,10.00){\line(0,1){8.00}}
\put(110.00,10.00){\line(0,1){8.00}}
\put(80.00,18.00){\line(1,0){30.00}}

\put(90.00,19.00){$\scriptstyle{{\rm color}\;q}$}
\put(70.50,5.00){$\scriptstyle{F_p\ell(u'')^*F_q}$}
\put(100.50,5.00){$\scriptstyle{F_q\ell(u'')F_p}$}
\end{picture}
\\
Now, the first one corresponds to $F_q\ell(u')^*F_pwF_p\ell(u')F_q$ for some word 
$w$, which produces $b_{p,q}(u')F_q$, and the second one corresponds to $F_p\ell(u'')^*F_quF_q\ell(u'')F_p$ 
for some word $u$, 
which produces $b_{q,p}(u'')F_p$, with $F_q$ and $F_p$ left to be matched 
by the neighboring diagonal matrices and forcing the nearest outer block to be 
colored by $q$ and $p$, respectively. We also use the fact that the covariances assigned to $u'$ and $u''$ are equal. It can be seen that these contributions correspond to the way covariances are assigned to blocks in the sense of Definitions 6.3 and 6.4. Therefore, our proof is completed.
\hfill $\blacksquare$\\

\section{Cumulants}

We would like to examine matricial circular operators from the point of view of cumulants. 
This leads to a family of scalar-valued cumulants labelled by $q$, whose definition 
is tailored to the matricial nature of our variables and is based on the cyclicity condition. 
A similar idea appeared 
in the context of $R$-cyclic families of matrices [9,10], except that we use a family of states 
rather than one trace.

\begin{Definition}
{\rm A family of multilinear functions $\kappa_{\pi}[\;.\,\;;q]$, where $q\in [r]$, of matricial 
variables 
$$
a_j=c_j\otimes e(p_j,q_j)\in M_{r}({\mathcal A}),
$$
where $1\leq j \leq m$, will be called {\it cyclically multiplicative} over 
the blocks $V$ of the partition $\pi\in \mathcal{NC}_{m}(a_1, \ldots, a_m)$ if 
$$
\kappa_{\pi}[a_1, \ldots , a_m;q]=\prod_{{\rm blocks}\;V}\kappa_{V}[a_1, \ldots , a_m;q_{V}],
$$
where 
$$
\kappa_{V}[a_1, \ldots , a_m;q_{V}]=\kappa_{s}(a_{i(1)}, \ldots , a_{i(s)};q_{i(s)})
$$
for the block $V=(i(1)<\ldots <i(s))$, i.e. $q_{V}$ is equal to the last matrix index 
which is associated with the block $V$, where $\{\kappa_{n}(.;q):n\geq 1, q\in [r]\}$ is a family of multilinear functions. If $\pi \notin \mathcal{NC}_{m}(a_1, \ldots ,a_m)$, then we set
$\kappa_{\pi}[a_1, \ldots , a_m;q]=0$ for any $q$.
}
\end{Definition}

\begin{Definition}
{\rm 
With the above notations, by the {\it cyclic cumulants} we shall understand 
the family of multilinear cyclically multiplicative functionals over the blocks of 
noncrossing partitions 
$$
\pi\rightarrow \kappa_{\pi}[\,.\,;q],
$$
defined by $r$ moment-cumulant formulas
$$
\Phi_q(a_1\ldots a_m)=\sum_{\pi\in \mathcal{NC}_{m}(a_1, \ldots ,a_m)}
\kappa_{\pi}[a_1, \ldots , a_m;q],
$$
where $q\in [r]$ and $\Phi_q=\varphi\otimes \psi_q$ are states on $M_{r}({\mathcal A})$ as defined before.}
\end{Definition}

Let us establish a relation between cyclic cumulants and operator-valued free cumulants.
For details on operator-valued cumulants, see [14].
\begin{enumerate}
\item
The $A$-{\it valued free cumulants} are $A$-functionals 
$\kappa_{m}:{\mathcal C}^{m}\rightarrow A$ 
defined recursively by the moment-cumulant formula
$$
E(a_1\ldots a_m)=\sum_{\pi\in \mathcal{NC}_{m}}\kappa_{\pi}[a_1, \ldots , a_m]
$$
where $\kappa_{\pi}[a_1, \ldots , a_m]$ are {\it partitioned $A$-valued cumulants}. 
\item
These are defined by setting $\kappa_{\pi}[a_1, \ldots , a_m]=\kappa_{m}(a_{1}, \ldots , a_{n})$ 
when $\pi=1_{m}$, which means that $\pi$  consists of one block, and 
requiring that otherwise they satisfy the recursion
$$
\kappa_{\pi}[a_1, \ldots , a_m]= 
\kappa_{\pi\setminus V}[a_1, \ldots ,a_j \kappa_{|V|}[a_{j+1}, \ldots , a_{k}],a_{k+1}, \ldots , a_m]
$$
whenever the interval $V=[j+1,k]$ is a block of $\pi$. 
\item
Recall that $\rho:{\mathcal C}^{n}\rightarrow A$ is an $A$-{\it functional} if it is 
a multilinear map such that 
$$
\rho[b_0a_1b_1,a_{2}b_{2}, \ldots , a_mb_m]=
b_{0}\rho[a_1, b_1a_{2},\ldots , b_{m-1}a_m]b_{m}
$$
for any $b_0,\ldots, b_m\in A$ and $a_{1}, \ldots ,a_{m}\in {\mathcal C}$. 
\item
We will use the above definitions in the case when ${\mathcal C}=M_{r}({\mathcal A})$ and 
the conditional expectation is of the form
$$
E(\sum_{p,q}a(p,q)\otimes e(p,q))
=\sum_{q}\varphi(c(q,q))\otimes e(q,q)
$$
where $A$ is again the algebra of $r\times r$ diagonal complex matrices.
\end{enumerate}

\begin{Proposition}
If $E:M_{r}({\mathcal A})\rightarrow A$ is as above, then 
$$
\kappa_{m}[a_1, \ldots , a_{m}]=\sum_{q=1}^{r}\kappa_{m}[a_1, \ldots , a_m; q]F_q
$$
for any $m\in \mathbb{N}$ and $a_{1}, \ldots , a_{m}\in M_{r}({\mathcal A})$.
\end{Proposition}
{\it Proof.}
It suffices to show that 
$$
\kappa_{\pi}[a_1, \ldots , a_m]=
\kappa_{\pi}[a_{1}, \ldots , a_{m};q]F_{q}
$$
for any $\pi\in {\mathcal NC}_{m}(a_{1}, \ldots , a_{m})$ with $q_{m}=q$, and otherwise 
$\kappa_{\pi}[a_1, \ldots , a_m]=0$.
It is clear that both sides vanish when $\pi$ is not adapted to $(a_1, \ldots, a_m)$.
Moreover, $\kappa_{\pi}[a_1]=E(a_{1})F_q=\kappa_{\pi}[a_1;1]F_{q}$ when $\pi=1_1$.
Suppose the equation holds for cumulants of order $<m$.
Now, if $\pi$ is adapted to $(a_1, \ldots, a_m)$ and $\pi\neq 1_{m}$, then 
it has a block $V=[j,k]$ such that $p_j=q_{k}=p$. Then
$$
\kappa_{\pi}[a_1, \ldots , a_m]= \kappa_{|V|}(a_{j+1}, \ldots , a_{k};q_k)
\kappa_{\pi\setminus V}[a_1, \ldots ,a_{j-1}, a_{k+1}, \ldots , a_{m}] 
$$
since $a_{j-1}F_{p_j}=a_{j-1}$ and $F_{q_k}a_{k+1}=a_{k+1}$ for $p_j=q_k$.
Namely, we can pull out the scalar-valued cyclic cumulant 
$\kappa_{V}(a_1, \ldots , a_m; q_k)$ and we do not need to keep 
$F_{q_{k}}$ inside. By induction, we obtain 
that $\kappa_{\pi\setminus V}[a_1, \ldots ,a_{j-1},a_{k+1}, \ldots,  a_m]$ is a product of 
lower-order cyclic cumulants corresponding to the remaining blocks of $\pi$ multiplied by 
$F_q$ since we can always choose $V$ so that $m\notin V$ and then clearly $q_m=q$. 
Therefore, we obtain exactly the expression for $\kappa_{\pi}[a_1, \ldots , a_m;q]$
in Definition 6.1 multiplied by $F_q$.
Finally, if $\pi$ is adapted to $(a_{1}, \ldots , a_{m})$ and $\pi=1_{m}$, then 
the moment-cumulant formula and induction gives 
$$
\kappa_{m}[a_{1}, \ldots , a_{m}]=\kappa_{m}(a_{1}, \ldots , a_m; q)F_q
$$
where $p_1=q_m=q$, and in consequence, the desired equation for partitioned cumulants. 
This completes the proof.
\hfill $\blacksquare$.\\

Treating all matricial circular operators as `independent' variables, we introduce the family 
of {\it cyclic $R$-transforms} of the arrays $\zeta:=(\zeta_{p,q})$  and 
$\zeta^*:=(\zeta_{p,q}^*)$. These are formal power series in $2r^2$ noncommuting indeterminates 
arranged in arrays $z:=(z_{p,q})$ and $z^*:=(z_{p,q}^{*})$ of the form
$$
R_{\zeta, \zeta^{*}}(z^{}, z^{*};q)=
\sum_{n=1}^{\infty}
\sum_{\stackrel
{p_1, \ldots , p_n}
{\scriptscriptstyle q_1, \ldots , q_n}
}
\sum_{\epsilon_1, \ldots, \epsilon_n}
\kappa_{n}(\zeta^{\epsilon_1}_{p_1,q_1}, 
\ldots , \zeta_{p_n,q_n}^{\epsilon_n};q)
z_{p_1,q_1}^{\epsilon_1}\ldots z_{p_n,q_n}^{\epsilon_n},
$$
where it is understood that $p_1, \ldots , p_n, q, q_1, \ldots q_n\in [r]$ and 
$\epsilon_1, \ldots , \epsilon_n\in \{1, *\}$. The theorem proved below shows that 
the cyclic $R$-transforms are the canonical transforms of $R$-transform type to describe 
circular systems.
 
\begin{Theorem}
If $\zeta:=(\zeta_{p,q})$ is the square array of matricial circular operators and 
$\zeta^*:=(\zeta^*_{p,q})$, then their cyclic $R$-transforms are of the form 
$$
R_{\zeta, \zeta^{*}}(z^{}, z^{*};q)=
\sum_{p=1}^{r}(b_{p,q}^{} z_{p,q}^{*}z_{p,q}^{}+ b_{p,q}^{}z_{q,p}^{}z_{q,p}^{*})
$$
for any $q\in [r]$, where 
$b_{p,q}=\kappa_{2}(\xi_{p,q}^{*},\xi_{p,q};q)=\kappa_{2}(\xi_{q,p},\xi_{q,p}^{*};q)$ 
for any $p,q$.
\end{Theorem}
{\it Proof.}
We present a proof based on operator-valued free cumulants,
Since we can identify $\zeta_{p,q}$ with $F_pgF_q$ and it is known that 
only second-order $A$-valued cumulants of $g$ do not vanish (as $g$ is a sum 
of two semicircular elements), we have
$$
\kappa_{2}(\zeta_{p,q}^*, \zeta_{p,q})=\kappa_{2}(\zeta_{q,p}, \zeta_{q,p}^*)=b_{p,q}F_q
$$
and thus, by Proposition 6.1, we get the desired formula for the cyclic $R$-transform.
\hfill $\blacksquare$\\

In the results given above, we dealt with cyclic $R$-transforms 
of square $r\times r$ arrays of matricial circular operators. However, one can also consider
pairs consisting of one fixed operator and its adjoint. In the off-diagonal case, 
we then obtain two non-trivial cyclic $R$-transforms of this pair and they take
a very simple form.

\begin{Corollary}
If $p,q\in [r]$ are fixed and $p\neq q$, then  
\begin{enumerate}
\item 
the non-trivial cyclic R-transforms of $h:=\zeta_{p,q}, h^{*}:=\zeta^{*}_{p,q}$ are
$$
R_{h, h^{*}}(z_1,z_2;q)=b_{p,q} z_2z_1\;\;\;and\;\;\;
R_{h, h^{*}}(z_1,z_2;p)=b_{q,p} z_1 z_2,
$$
\item
the non-trivial cyclic R-transform of $c:=\zeta_{q,q}, c^{*}=\zeta_{q,q}^{*}$ is 
$$
R_{c, c^{*}}(z_1,z_2;q)=b_{q,q}(z_1z_2+z_2z_1),
$$
\item
the non-trivial cyclic R-transforms of $\eta:=\eta_{p,q}, \eta^{*}:=\eta^{*}_{p,q}$ are
\begin{eqnarray*}
R_{\eta, \eta^{*}}(z_1,z_2;q)&=&b_{p,q}(z_1z_2+z_2z_1)\\
R_{\eta, \eta^{*}}(z_1,z_2;p)&=&b_{q,p}(z_1z_2+z_2z_1).
\end{eqnarray*}
\end{enumerate}
\end{Corollary}
{\it Proof.}
The formulas for these cyclic $R$-transforms are obtained from Theorem 6.1 
by taking suitable subarrays.
\hfill $\blacksquare$\\

\begin{Remark}
{\rm
It seems natural that the cyclic $R$-transforms of $\eta_{p,q},\eta_{p,q}^*$ in the states $\Phi_q$
and $\Phi_p$ coincide with the $R$-transforms of circular operators since the *-distributions of 
$\eta_{p,q}$ in these states are circular. However, cyclic $R$-transforms are defined 
differently than the usual $R$-transforms and thus we prefer to treat this property 
as a decomposition of the $R$-transform of circular operators in terms of cyclic $R$-transforms 
of matricial circular operators rather than as an obvious fact. A more general decomposition property
of this type is formulated below.}
\end{Remark}

\begin{Corollary}
If $b_{p,q}=d_{p}$ for any $p,q$, then $c:=\sum_{p,q}\zeta_{p,q}$ is circular with respect to 
the convex linear combination $\Phi=\sum_{q=1}^{r}d_{q}\Phi_q$ and the corresponding 
$R$-transform takes the form
$$
R_{c,c^{*}}(z_1,z_2)=\sum_{q=1}^{r}d_q R_{\zeta,\zeta^*}(z,z^*;q),
$$
where $R_{\zeta,\zeta^*}(z,z^*;q)$ is the cyclic $R$-transform of Theorem 8.1 in which 
all entries of arrays $z$ and $z^{*}$ are identified with $z_1$ and $z_2$, respectively. 
\end{Corollary}
{\it Proof.}
If we assume that all block variances in [7, Theorem 9.1] are equal to one, the moments of the matrix $Y(n)$ 
converge under the normalized trace to those of the standard circular operator 
$c:=\sum_{p,q}\zeta_{p,q}$ under $\Phi$. We have
$$
\kappa_{2}(c^*,c)=\Phi(c^*c)\;\;\;{\rm and}\;\;\;\kappa_{2}(\zeta_{p,q}^{*}, \zeta_{p,q}^{};q)=\Phi_{q}(\zeta_{p,q}^{*}\zeta_{p,q}^{}),
$$
where we denote by $\kappa_{n}$ the (scalar-valued) free cumulants with respect to $\Phi$. 
Since
$$
\Phi(c^*c)=\sum_{p,q=1}^{r}d_{q}\Phi_{q}(\zeta_{p,q}^{*}\zeta_{p,q}^{}),
$$
we obtain 
$$
\kappa_{2}(c^*,c)=\sum_{q=1}^{r}d_{q}\kappa_{2}(c^*,c;q).
$$
A similar decomposition holds if we interchange $c$ and $c^*$.
Free cumulants of $c,c^*$ of higher orders vanish since $c$ is circular.
The cyclic cumulants of higher orders of $c,c^*$ also vanish by Proposition 6.1. 
Therefore, the $R$-transform for $c,c^*$ has the desired decomposition in terms of their 
cyclic $R$-transforms which can be identified with 
$R_{\zeta^{}, \zeta^*}(z,z^*;q)$ for $q\in [r]$, where $z_{p,q}=z_1$ and $z_{p,q}^{*}=z_2$
for any $p,q$. Of course, both sides of this equation are equal to $z_1z_2+z_2z_1$.
\hfill $\blacksquare$\\

\section{Moment series for matricial circular systems}

The matricial moment series for the arrays $\zeta=(\zeta_{p,q})$ and $\zeta^*=(\zeta_{p,q}^{*})$ is the family of formal series of the form
$$
M_q(z,z^*)=\sum_{m=1}^{\infty}
\sum_{\stackrel{p_1,q_1, \ldots , p_m,q_m}
{\scriptscriptstyle \epsilon_1, \ldots , \epsilon_m}
}
\Phi_q(\zeta_{p_1,q_1}^{\epsilon_1}\ldots \zeta_{p_m,q_m}^{\epsilon_m})z_{p_1,q_1}^{\epsilon_1}\ldots z_{p_m,q_m}^{\epsilon_m}
$$
where $q\in [r]$ and $z=(z_{p,q})$ and $z^*=(z^*_{p,q})$ are the arrays of 
$2r^2$ noncommuting indeterminates. The moment series for a family of arrays labelled by 
$u$ is defined in a similar way. 
All results of this Section can be generalized accordingly.

\begin{Proposition}
The moment series for the arrays $\zeta, \zeta^*$ satisfy the relations
$$
M_{q}^{}=1+\sum_{p=1}^{r}b_{p,q}(z_{p,q}^{*}M_{p}^{}z_{p,q}^{}M_{q}^{}+
z_{q,p}^{}M_{p}^{}z_{q,p}^{*}M_{q}^{})
$$
for any $q\in [r]$, where $M_{s}=M_{s}(z,z^*)$ for any $s\in [r]$. 
\end{Proposition}
{\it Proof.}
We can express each moment $\Phi_q(\zeta_{p_1,q_1}^{\epsilon_1}\ldots \zeta_{p_m,q_m}^{\epsilon_m})$
in terms of cyclic cumulants according to the formula
$$
\Phi_{q}(\zeta_{p_1,q_1}^{\epsilon_1}\ldots \zeta_{p_m,q_m}^{\epsilon_m})
=
\sum_{\pi\in \mathcal{NC}_{m}^{2}(\zeta_{p_1,q_1}^{\epsilon_1}, \ldots, \,\zeta_{p_m,q_m}^{\epsilon_m})}
\kappa_{\pi}[\zeta_{p_1,q_1}^{\epsilon_1}, \ldots, \zeta_{p_m,q_m}^{\epsilon_m};q],
$$
where each partitioned cyclic cumulant $\kappa_{\pi}[.;q]$  is a product of cyclic cumulants over the 2-blocks of $\pi$.
In order to find a typical relation between moment series, it now suffices to single out in each $\pi$ 
the 2-block containing $1$, say $\{1,k\}$. The cyclic cumulant corresponding to this block is $b_{p,q}$ for some $p$
and the corresponding matrix units must be $e(q,p)$ and $e(p,q)$, respectively.    
All cyclic cumulants corresponding to the inner blocks of this block are collected in  
the series $M_{p}(z,z^*)$ since they correspond to the mixed *-moments of variables lying in between 
$\zeta_{p_1,q_1}^{\epsilon_1}$ and $\zeta_{p_k,q_k}^{\epsilon_k}$ computed `under condition $p$' since the corresponding product 
of matrix units must be $e(p,p)$. In turn, all cyclic cumulants correponding to the blocks involving numbers greater than $q$ are collected in the series
$M_{q}(z,z^*)$ since they correspond to the mixed *-moments of $\zeta_{p_{k+1},q_{k+1}}^{\epsilon_{k+1}}, \ldots, \zeta_{p_m,q_m}^{\epsilon_m}$ 
and the corresponding product of matrix units must be $e(q,q)$. Of course, $z_{p,q}$ is assigned to $\zeta_{q,p}^*$ and $z_{p,q}$ is assigned 
to its adjoint $\zeta_{p,q}^{*}$ and, since our indeterminates do not commute, we have 
to place the series $M_{p}(z,z^*)$ in between them if $\{1,k\}$ is associated with the pairing $(\zeta_{p,q}^*, \zeta_{p,q})$.
In turn, the series $M_{p}(z,z^*)$ is placed in between $z_{q,p}$ and $z_{p,q}^*$ if $\{1,k\}$ is associated with the 
pairing $(\zeta_{q,p}, \zeta_{q,p}^{*})$. This gives the formula
$$
M_{q}=1+\sum_{p=1}^{r}b_{p,q}^{}(z_{q,p}^{*}M_{p}z_{p,q}^{}M_{q}+
z_{q,p}^{}M_{p}z_{p,q}^{*}M_{q})
$$
where $M_{s}=M_{s}(z,z^*)$ for any $s$.
\hfill $\blacksquare$\\

One can express the formula of Proposition 7.1 in terms of the cyclic $R$-transforms of $\zeta, \zeta^*$. 
For that purpose, it is convenient to replace the array $z^*$ by its transpose $w=(w_{p,q})$, where 
$w_{p,q}=z_{q,p}^*$ for any $p,q\in [r]$. Moreover, let
\begin{eqnarray*}
M(z,w):&=&{\rm diag}(M_1, \ldots , M_r)\\
R(z,w):&=&{\rm diag}(R_1, \ldots , R_r),
\end{eqnarray*}
where $M_q=M_q(z,w)$ and $R_q=R_{\zeta, \zeta^*}(z,w;q)$
for any $q\in [r]$ and we abuse the notation in the sense that we replace each $z_{p,q}^*$ by $w_{q,p}$ and 
we still write $w$ as the second argument of each $M_q$ and $R_q$.

\begin{Corollary}
With the above notations, let
$$
R_{q}(z,w)=\sum_{p=1}^{r}b_{p,q}(w_{q,p}z_{p,q}+
z_{q,p}w_{p,q})
$$
for any $q\in [r]$. Then, it holds that 
$$
M(z,w)=I+R(zM(z,w), wM(z,w)).
$$ 
where $zM(z,w)$ and $wM(z,w)$ are products of arrays.
\end{Corollary}
{\it Proof.}
We can substitute $z_{s,t}\rightarrow z_{s,t}M_{t}=(zM)_{s,t}$ and $w_{s,t}\rightarrow w_{s,t}M_{t}=(wM)_{s,t}$ 
for $s,t\in \{p,q\}$, which gives the formula of Proposition 7.1. This completes the proof.
\hfill $\blacksquare$\\

From the formula of Proposition 7.1 we can obtain a concise formula for the diagonal matrix of moment series 
for $\zeta=\sum_{p,q}\zeta_{p,q}$ in the states $\Psi_q$. Let $\mathcal{D}:M_{n}({\mathcal A})\rightarrow M_{n}({\mathcal A})$ 
be the linear mapping given by 
$$
\mathcal{D}(A)={\rm diag}(A_1, \ldots , A_r)
$$ 
for any $A\in M_{r}({\mathcal A})$, where
$$
A_q=\sum_{p=1}^{r}a_{p,q}
$$
for any $q$. This formula will be similar to that derived in [7] for the moment series
for the sums of matricial semicircular operators $\omega=\sum_{p,q}\omega_{p,q}$. 

\begin{Corollary}
Let $M={\rm diag}(M_1, \ldots , M_r)$, where each $M_q$ is given by  
$$
M_q(z,z^*)=\sum_{m=1}^{\infty}\sum_{\epsilon_1, \ldots , \epsilon_m} \Phi_q(\zeta^{\epsilon_1}\ldots \zeta^{\epsilon_{m}})
z^{\epsilon_1}\ldots z^{\epsilon_{m}}
$$
for the sum $\zeta=\sum_{p,q}\zeta_{p,q}$. Then $M$ satisfies the functional equation
$$
M=I+\mathcal{D}(z^*MBzM+zMBz^*M),
$$
where $B=(b_{p,q})$ and $z,z^*$ are two noncommuting variables.
\end{Corollary}
{\it Proof.}
We set $z_{s,t}=z$ and $z_{s,t}^{*}=z^{*}$ for any $s,t$ in the 
formulas of Proposition 7.1, which gives
$$
M_{q}=1+\sum_{p=1}^{r}(z^{*}M_{p}b_{p,q}zM_{q}+
zM_{p}b_{p,q}z^{*}M_{q})
$$
for any $q$, where we put $b_{p,q}$ between $M_p$ and $M_q$ to anticipate the product of matrices $M_p, B, M_q$. Now, 
it is easy to see that $z^*M_{p}b_{p,q}zM_{q}$ is the $(p,q)$-th entry of the matrix $z^*MBzM$, where
$z$ and $z^*$ are to be understood as diagonal matrices with $z$ and $z^*$ on the diagonal, respectively.
Similarly, $zM_pb_{p,q}z^*M_q$ is the $(p,q)$-th entry of $zMBz^*M$. In turn, summation over $p$ for fixed $q$ 
corresponds to the mapping $\mathcal{D}$, which completes the proof. 
\hfill $\blacksquare$\\

The formulas of Proposition 7.1 for the moment series for arrays of $R$-circular operators as well as those of Corollary 7.2 
for the moment series for their sums are equivalent to the definition of their cyclic cumulants (which are recursions for moments) 
and thus it is not surprising that their solutions can be found in the form of recursions analogous to those for Catalan numbers. 
For simplicity, we restrict our attention to the case discussed in Corollary 7.1.

\begin{Proposition}
The matrix-valued moment series $M$ of Corollary 7.2 takes the form 
$$
M(z,z^*)=\sum_{m=0}^{\infty}\sum_{\epsilon_1, \ldots , \epsilon_{2m}}C_{m}(\epsilon_1, \ldots, \epsilon_{2m})\,
z^{\epsilon_1}\ldots z^{\epsilon_{2m}},
$$
where the coefficients are diagonal $r\times r$ matrices satisfying the recurrence
$$
C_{m}(\epsilon_1, \ldots, \epsilon_m)=\sum_{j+k=m-1}\mathcal{D}(C_{j}(\epsilon_2, \ldots , \epsilon_{2j-1})BC_{k}(\epsilon_{2j+1}, \ldots , \epsilon_{2m}))
$$
for any $m\in {\mathbb N}$, with $C_{0}=I$, where we set $C_{m}(\epsilon_1, \ldots, \epsilon_m)=0$ whenever
$(\epsilon_1, \ldots, \epsilon_m)$ is not associated with a noncrossing pairing consisting of $(1,*)$ and $(*,1)$. 
This recurrence has a unique solution.
\end{Proposition}
{\it Proof.}
It suffices to substitute $M(z,z^*)$ in the above form to the equation for $M$ of Corollary 7.2 and observe that if
$(\epsilon_1, \epsilon_{2j})=(*,1)$, then $C_j(\epsilon_2, \ldots , \epsilon_{2j-1})$
comes from the first $M$ in the product $z^*MBzM$ and $C_{k}(\epsilon_{2j-1}, \ldots , \epsilon_{2m})$ 
comes from the second $M$ in that product. All terms, for which $(\epsilon_1, \epsilon_{2j})=(1,*)$,
correspond to the product $zMBz^*M$ in a similar way. This gives the desired recurrence for the matrix-valued coefficients.
It is obvious that this recurrence has a unique solution.
\hfill$\blacksquare$ \\

\begin{Corollary}
If $b_{p,q}=d_p$ for any $p,q\in [r]$, then the 
moment series $M_{0}=M_{0}(z,z^*)$ for the circular operator $\zeta=\sum_{p,q}\zeta_{p,q}$ 
associated with $\Phi=\sum_{q=1}^{r}d_{q}\Phi_q$ satisfies the equation
$$
M_0=1+z^*M_0zM_0+zM_0z^*M_0
$$
\end{Corollary}
{\it Proof.}
It holds that 
$$
M_0(z,z^*)=\sum_{m=1}^{\infty}\sum_{\epsilon_1, \ldots , \epsilon_m} 
\Phi(\zeta^{\epsilon_1}\ldots \zeta^{\epsilon_{m}})
z^{\epsilon_1}\ldots z^{\epsilon_{m}}=\sum_{q=1}^{r}d_{q}M_q(z,z^*)
$$
and it suffices to take the sum over $q$ of the equation in the proof 
of Corollary 7.2 and set $b_{p,q}=d_p$ to obtain the desired equation.
\hfill $\blacksquare$\\

\section{Circular free Meixner distributions}

Using matricial circular systems, we will first define operators
whose distributions play the role of non-Hermitian counterparts of Kesten distributions. 
We call them `circular Kesten distributions'. A natural generalization of circular Kesten distributions
leads to `circular free Meixner distributions' which play the role of non-Hermitian counterparts 
of free Meixner distributions.

Using Theorem 4.1, we will show that circular Kesten distributions are 
limit *-distributions of non-Hermitian GRM (with i.b.i.d. entries) of the block 
form 
$$
Y(n)=\left(
\begin{array}{rr}
A(n)& B(n)\\
C(n)& D(n)
\end{array}
\right),
$$
under the partial trace $\tau_1(n)$ as $n\rightarrow \infty$. These matrices 
satisfy the assumptions of Theorem 4.1 with asymptotic dimensions $d_{1}=0$ and $d_2=1$. 
This corresponds to the situation in which $(A(n))$ is {\it evanescent}, $(D(n))$ is {\it balanced}, 
whereas $(B(n))$ and $(C(n))$ are {\it unbalanced} for any $u$. This is the non-Hermitian counterpart 
of the random matrix model for Kesten laws constructed recently in [8]. 

When discussing *-distributions of one such sequence of matrices (or, *-distributions of their blocks) we omit $u\in \mathpzc{U}$ in our notations. More generally, we can take a family of non-Hermitian GRM of the above type, 
$$
\{Y(u,n):u\in \mathpzc{U}, n\in \mathbb{N}\},
$$ 
and their blocks. Then their limit joint *-distribution 
is expressed in terms of the whole matricial circular system 
of Definition 3.1 instead of one array $(\zeta_{p,q})$.

\begin{Definition}
{\rm 
If $(\zeta_{p,q})_{1\leq p,q\leq 2}$ is an array of matricial circular operators
such that $\zeta_{2,1}$ and $\zeta_{1,2}^*$ have covariances $\beta_1 $ and $0$ 
in the states $\Psi_1$ and $\Psi_2$, respectively, whereas $\zeta_{2,2}$ has covariances 
$0$ and $\beta_2$ in the states $\Psi_1,$ and $\Psi_2$, respectively, 
where $\beta_1, \beta_2>0$, then the operator of the form
$$
\zeta=\zeta_{2,1}+\zeta_{1,2}+\zeta_{2,2},
$$ 
will be called a {\it circular Kesten operator}. Its *-distribution w.r.t. $\Psi_1$
will be called the {\it circular Kesten distribution}
associated with $(\beta_1, \beta_2)$. If $\beta_1=\beta_2=1$, it reduces to 
the circular operator with the standard circular distribution w.r.t. $\Psi_1$. 
}
\end{Definition}

\begin{Remark}
{\rm We justify the above terminology by comparing circular Kesten operators
with the self-adjoint operators having Kesten distributions studied in [8].
\begin{enumerate}
\item
If $\omega_{2,1},\omega_{2,2}$ have variances $\beta_1, \beta_2$ in the 
states $\Psi_1, \Psi_2$, respectively, then the operator of the form
$$
\omega=\omega_{2,1}+\omega_{2,2},
$$
has the Kesten distribution associated with $(\beta_1, \beta_2)$.
In fact, it corresponds to the sequence 
of Jacobi coefficients $(\beta_1, \beta_2, \beta_2, \ldots )$.
Equivalently, we can write
\begin{eqnarray*}
\omega
&=&
\wp_{2,1}+\wp_{2,1}^*+\wp_{2,2}+\wp_{2,2}^*\\
&=&
\ell_1\otimes e(2,1)+\ell_1^*\otimes e(1,2)+(\ell_2+\ell_2^*)\otimes e(2,2)
\end{eqnarray*}
where $\ell_1, \ell_2$ are free under $\varphi$ and have covariances $\beta_1, \beta_2$, respectively.
\item
Using the fact that $\wp_{1,2}(u')=\wp_{1,2}(u'')=0$ (by the assumption on zero covariances in $\Psi_2$), circular Kesten operators take the form
\begin{eqnarray*}
\zeta
&=& 
\wp_{2,1}(u')+\wp_{2,1}^*(u'')+\wp_{2,2}(u')+\wp_{2,2}^*(u'')\\
&=&
\ell_1\otimes e(2,1)+\ell_3^*\otimes e(1,2)+(\ell_2+\ell_4^*)\otimes e(2,2),
\end{eqnarray*}
where $u'\neq u''$ and $\ell_1,\ell_3$ have both covariance $\beta_1$ and 
$\ell_2, \ell_4$ have both covariance $\beta_2$ in the state $\varphi$ 
(Proposition 3.1 is used). This expression 
shows that each $\zeta$ can be viewed as non-self-adjoint analog of 
some $\omega$.
\end{enumerate}}
\end{Remark}

It turns out that the combinatorics of circular Kesten distributions 
generalizes that of circular distributions inasmuch as the 
combinatorics of Kesten distributions generalizes 
that of semicircular distributions.

By $\mathcal{NC}^{2}_{m}(\epsilon_1, \ldots ,\epsilon_m)$ 
we denote the set of noncrossing pair-partitions of $\{1, \ldots , m\}$ 
which are adapted to $(\epsilon_1, \ldots , \epsilon_m)$, that is 
$V=\{i,j\}\in \pi$ if and only if $\epsilon_i\neq \epsilon_j$.
Moreover, $\mathpzc{d}(V)$ stands for the depth of block $V$ and we set $\beta_{{\mathpzc d}}=\beta_2$ also for ${\mathpzc d}> 2$. 

\begin{Proposition}
The *-distribution of the circular Kesten operator
$\zeta$ associated with $(\beta_1, \beta_2)$ in the state $\Psi_1$ is given by mixed *-moments 
$$
\Psi_1(\zeta^{\epsilon_1}\ldots \zeta^{\epsilon_m})=
\sum_{\pi\in \mathcal{NC}^{2}_{m}(\epsilon_1, \ldots , \epsilon_m)}
\beta(\pi)
$$
where
$$
\beta(\pi)=\prod_{{\rm blocks}\;V}\beta_{\mathpzc{d}(V)}
$$
and the product is taken over all blocks of $\pi$. 
\end{Proposition}
{\it Proof.}
It is clear that the given *-moment is a sum of mixed *-moments of the form 
$$
\Psi_1(\zeta_{p_1,q_1}^{\epsilon_1}\ldots \zeta_{p_m,q_m}^{\epsilon_m})
$$
where $(p_1,q_1), \ldots , (p_m,q_m)\in \{(1,2),(2,1), (2,2)\}$.
We then employ the combinatorial formula of Lemma 5.1. Since we use the state $\Psi_1$, we consider noncrossing 
colored pair-partitions with the imaginary block colored by $1$. Now, the covariances of
$\wp_{p,q}(u')$ and $\wp_{p,q}(u'')$ in the state $\Psi_q$
are both equal to $\beta_{q}$ for $q\in \{1,2\}$ (if we write $\zeta=\zeta(u)$,
then $u'$ and $u''$ have the same meaning as in Definition 3.1). 
Thus, to each block $\{k,l\}$ of
a colored pair-partition $\pi$ which is adapted to the tuple $(\epsilon_1, \ldots , \epsilon_m)$
there corresponds the pair $(\zeta_{p_k,q_k}^{\epsilon_k}, \zeta_{p_l,q_l}^{\epsilon_l})$, where
$(p_k,q_k)=(p_l,q_l)$ and $\epsilon_k\neq \epsilon_l$, which 
contributes $\beta_1=b_{2,1}(u)$ when $\{k,l\}$ is of depth 1 and 
$\beta_2=b_{2,2}(u)$ when $\{k,l\}$ is of depth $\geq 2$. There is a bijection between 
non-vanishing mixed *-moments of this type and colored noncrossing pair-partitions adapted to
$(\epsilon_1, \ldots , \epsilon_m)$, which gives the desired formula. 
\hfill $\blacksquare$

\begin{Theorem}
Let $\beta_1=v_{1,2}=v_{2,1}>0$ and $\beta_2=v_{2,2}>0$ be the 
covariances of variables in blocks $B(n), C(n)$ and $D(n)$, respectively, 
in a Gaussian random matrix $Y(n)$ with i.b.i.d. entries. 
If the asymptotic dimensions are $d_1=0$ and $d_2=1$, then
$$
\lim_{n\rightarrow \infty}\tau_1(n)\left(Y(n)^{\epsilon_1}\ldots Y(n)^{\epsilon_m}\right)=
\Psi_1(\zeta^{\epsilon_1}\ldots \zeta^{\epsilon_m})
$$
for any $m\in {\mathbb N}$ and any $\epsilon_1, \ldots , \epsilon_m\in \{1,*\}$,
where $\zeta$ is the circular Kesten operator associated with $(\beta_1, \beta_2)$.
\end{Theorem}
{\it Proof.}
By the assumption on asymptotic dimensions, the sequence of blocks $(A(n))$ is evanescent and 
thus $Y(n)\rightarrow \zeta$ as $n\rightarrow \infty$ by Theorem 4.1, 
where $\zeta$ is the circular Kesten operator associated with $(\beta_1, \beta_2)$, which proves our assertion.
\hfill $\blacksquare$\\

\begin{Remark}
{\rm One can generalize the above results to obtain
non-Hermitian counterparts of free Meixner laws. 
\begin{enumerate}
\item
Recall that in the case of free Meixner laws corresponding
to $\beta_1>0, \beta_2>0, \alpha_1, \alpha_2\in {\mathbb R}$, 
we used operators
$$
\gamma'=(\alpha_{2}-\alpha_1)(\beta_{1}^{-1}n_{2,1}+\beta_{2}^{-1}n_{2,2})+\alpha_{1}
$$
to construct operators $\xi'=\omega+\gamma'$ whose distributions 
in the state $\Psi_1$ were free Meixner laws [9, Theorem 3.1]. 
\item
Let us first suppose that $\beta_1>0$ and $\beta_2>0$. Introduce operators
$$
\gamma=(\alpha_{2}-\alpha_1)(\beta_{1}^{-1}(n_{2,1}(u')+n_{2,1}(u''))
+\beta_{2}^{-1}(n_{2,2}(u')+n_{2,2}(u'')))+\alpha_{1},
$$
where operators
$$
n_{p,q}(s)=\wp_{p,q}(s)\wp_{p,q}^{*}(s)
$$
play the role of rescaled `number operators' for any $p,q,s$
and $\alpha_1, \alpha_2\in {\mathbb C}$. Note that
$\gamma$ multiplies the vacuum vector by $\alpha_1$ and the remaining basis vectors 
of ${\mathcal M}$ by $\alpha_2$. 
\item
If $\beta_1=0$, then we set $\alpha_2=\beta_2=0$ and 
$\gamma=\alpha_1$.
\item
If $\beta_1>0$ and $\beta_2=0$, then we set 
$\alpha_2=0$ and
$$
\gamma=\alpha_1(1-\beta_{1}^{-1}n_{2,1}(u')-\beta_1^{-1}n_{2,1}(u'')).
$$
\end{enumerate}}
\end{Remark}

\begin{Definition}
{\rm 
Under the assumptions of Definition 8.1 and Remark 8.2, 
the operator of the form
$$
\xi=\zeta+\gamma
$$ 
will be called a {\it circular free Meixner operator}.
Its *-distribution in the state $\Psi_1$
will be called the {\it circular free Meixner distribution}
associated with $(\beta_1, \beta_2, \alpha_1, \alpha_2)$.
}

\end{Definition}
In order to write a combinatorial formula for the *-moments of circular free Meixner operators, 
we shall need the set $\mathcal{NC}^{1,2}_{m}(\epsilon_1, \ldots ,\epsilon_m)$ 
of noncrossing partitions of $\{1, \ldots , m\}$ consisting of singletons and 
pairs which are adapted to $(\epsilon_1, \ldots , \epsilon_m)$. Here, the notion of
adaptedness affects only the pairs of $\pi$, namely $V=\{i,j\}\in \pi$ if and only if 
$\epsilon_i\neq \epsilon_j$ and there are no restrictions on 
$\epsilon_i$ if $V=\{i\}$. Moreover, we set $\alpha_{{\mathpzc d}}=\alpha_2$ and 
$\beta_{{\mathpzc d}}=\beta_2$ also for ${\mathpzc d}> 2$. 

\begin{Proposition}
The *-distribution of the circular free Meixner operator
$\zeta$ associated with $(\beta_1, \beta_2, \alpha_1, \alpha_2)$ 
in the state $\Psi_1$ is given by mixed *-moments 
$$
\Psi_1(\xi^{\epsilon_1}\ldots \xi^{\epsilon_m})=
\sum_{\pi\in \mathcal{NC}^{1,2}_{m}(\epsilon_1, \ldots , \epsilon_m)}
\alpha(\pi)\beta(\pi)
$$
where 
$$
\alpha(\pi)=\prod_{{\rm singletons}\;V}\alpha_{\mathpzc{d}(V)}^{\epsilon(V)}
\;\;\;
and
\;\;\;
\beta(\pi)=\prod_{{\rm pairs}\;V}\beta_{\mathpzc{d}(V)}
$$
with products taken over sets of all singletons and pairs of $\pi$
and $\alpha^{*}=\overline{\alpha}$, where $\epsilon(V)=\epsilon_k$ if $V=\{k\}$.
\end{Proposition}
{\it Proof.}
The proof is similar to that of Proposition 8.1 (see also that for free Meixner laws
[8, Theorem 3.1]) and is omitted.
\hfill $\blacksquare$

\begin{Theorem}
Under the assumptions of Theorem 8.1 and for all $(\beta_1, \beta_2, \alpha_1, \alpha_2)$
considered above,
$$
\lim_{n\rightarrow \infty}\tau_1(n)\left(M(n)^{\epsilon_1}\ldots M(n)^{\epsilon_m}\right)
=\Psi_1(\xi^{\epsilon_1}\ldots \xi^{\epsilon_m}),
$$
where 
$$
M(n)=Y(n)+\alpha_1I_1(n)+\alpha_2I_2(n)
$$
for any $n\in {\mathbb N}$, where $I(n)=I_1(n)+I_2(n)$ is the decomposition of the 
$n\times n$ unit matrix induced by the partition $[n]=N_1\cup N_2$.
\end{Theorem}
{\it Proof.}
The proof is similar to that for free Meixner laws [8, Theorem 4.1] and is omitted.
\hfill $\blacksquare$\\ 

\begin{Remark}
{\rm The asymptotics of families of matrices of the 
types considered above can be obtained by a straightforward modification. 
It suffices to label each matrix and the corresponding operator by the same 
$u\in \mathpzc{U}$. 
\begin{enumerate}
\item
In Theorem 8.1, we obtain 
$$
Y(u,n)^{\epsilon}\rightarrow \zeta(u)^{\epsilon}
$$
under $\tau_{1}(n)$ as $n\rightarrow \infty$, where 
$\zeta(u)=\zeta_{2,1}(u)+\zeta_{1,2}(u)+\zeta_{2,2}(u)$ for any $u$ and $\epsilon$.
\item
In Theorem 8.2, we obtain  
$$
M(u,n)^{\epsilon}\rightarrow \xi(u)^{\epsilon}
$$
under $\tau_{1}(n)$ as $n\rightarrow \infty$, where 
$\xi(u)=\zeta(u)+\gamma(u)$ for any $u$ and $\epsilon$.
\item
In Propositions 8.1 and 8.2, we have to strenghten 
the notion of adaptedness for the considered 
noncrossing partitions of $\{1, \ldots , m\}$. 
Namely, they need to be adapted not only to 
$(\epsilon_1, \ldots , \epsilon_m)$, but also to 
$(u_1, \ldots , u_m)$, which means that $u_i=u_j$ whenever 
$\{i,j\}$ is a block of $\pi$ (singletons are not affected by this 
additional condition).
\item
These results can also be described, using the framework of operator-valued 
free probability, with convergence under partial traces replaced by 
convergence under conditional expectation. Moreover, circular Kesten operators
can be identified with $A$-{\it valued circular Kesten elements}
$$
\zeta^{A}(u)=F_1\ell(u')F_2+F_2\ell(u'')^*F_1+F_2g(u)F_2,
$$
where each $g(u)$ is an $A$-valued circular element. 
\item
Similarly, we can identify circular free Meixner operators with $A$-{\it valued circular free Meixner elements}
$$
\xi^{A}(u)=\zeta^{A}(u)+\gamma^{A}(u)
$$
where the expression for $\gamma^{A}(u)$ is similar to that for $\gamma(u)$, except that number operators  
$n_{p,q}(s)$ are replaced by $F_p\ell(s)F_q\ell(s)^*F_p$, respectively.
\end{enumerate}}
\end{Remark}

\noindent\\[10pt]
{\bf Acknowledgement}\\
I would like to express my thanks to the Reviewer for inspiring remarks which 
were very helpul in the preparation of the revised version of the paper.

\end{document}